\documentclass[12pt]{amsart}
			 
\usepackage{amsfonts}
\usepackage{amssymb}
\usepackage{amscd, graphics}
\usepackage{epsfig}
\input xy
\xyoption{all}

\newcommand{\marginlabel}[1]%
  {\mbox{}\marginpar{\raggedleft\hspace{0pt}\bfseries\sf#1}}

\def\ZZ{{\mathbb Z}}

\def\QQ{{\mathbb Q}}
\def\PP{{\textbf P}}
\def\OO{{\mathcal O}}

\def\cF{\mathcal{F}}

\def\cO{\mathcal{O}}
\def\cI{\mathcal{I}}
\def\I{\mathcal{I}}
\def\cM{\mathcal{M}}
\def\cZ{\mathcal{Z}}
\def\cU{\mathcal{U}}
\def\cQ{\mathcal{Q}}
\def\cE{\mathcal{E}}

\def\Mg{\cM_{g}}
\def\Mgn{\cM_{g,n}}
\def\Mg1{\cM_{g,g+1}}
\def\Pici{{\rm Pic}^{g-2i-1}(C)}

\def\iso{\simeq}

\def\tensor{\otimes}

\DeclareMathOperator{\Pic}{Pic}

\theoremstyle{plain}
\newtheorem*{introtheorem}{Theorem}   
\newtheorem*{introcorollary}{Corollary}
\newtheorem{theorem}{Theorem}[section]
\newtheorem{proposition}[theorem]{Proposition}
\newtheorem{corollary}[theorem]{Corollary}
\newtheorem{lemma}[theorem]{Lemma}

\theoremstyle{definition}
\newtheorem{definition}[theorem]{Definition}
\newtheorem{remark}[theorem]{Remark}
\newtheorem{example}[theorem]{Example}

\newtheorem{conjecture/question}[theorem]{Conjecture/Question}

\theoremstyle{remark}

\begin{document}

\title{Divisors on $\Mg1$ and the Minimal Resolution Conjecture
for points on canonical curves}

\author[G. Farkas]{Gavril Farkas}
\address{Department of Mathematics, University of Michigan, 
525 East University, Ann Arbor, MI, 48109-1109}
\email{{\tt gfarkas@umich.edu}}

\author[M. Musta\c{t}\v{a}]{Mircea~Musta\c{t}\v{a}}
\address{Department of Mathematics, University of California,
Berkeley, CA, 94720 and Institute of Mathematics of
the Romanian Academy}
\email{{\tt mustata@math.berkeley.edu}}

\author[M. Popa]{Mihnea Popa} 
\address{Department of Mathematics, Harvard University, 
One Oxford Street, Cambridge, MA 02138}
\email{{\tt mpopa@math.harvard.edu}}

\maketitle
\markboth{G. FARKAS, M. MUSTA\c T\v A and M. POPA}
{\bf THE MINIMAL RESOLUTION CONJECTURE}

\section*{\bf Introduction}

The Minimal Resolution Conjecture for points in projective space has 
attracted considerable attention in recent years, starting with the 
original \cite{lorenzini1}, \cite{lorenzini2} and continuing most notably 
with \cite{gaeta}, \cite{BG}, \cite{walter},
\cite{simpson}, \cite{popescu}, \cite{epsw}. The purpose
of this paper 
is to explain how a completely analogous problem can be formulated for 
sets of points on arbitrary varieties embedded in projective space, and 
then study in detail the case of curves. Similarly to 
the well-known analysis of syzygies of curves carried out by Green and 
Lazarsfeld (\cite{gl1}, \cite{gl2}, \cite{gl3}), we divide our work into
a study of resolutions of points on canonical curves and on curves of 
large degree. The central result of the paper states that the Minimal 
Resolution Conjecture is true on any canonical curve. In contrast, it 
always fails for curves embedded with large degree, although a weaker 
result, called the Ideal Generation Conjecture, holds also in this case.
These results turn out to have surprisingly deep connections with the 
geometry of difference varieties in Jacobians, special divisors on 
moduli spaces of curves with marked points, and moduli spaces of stable
bundles.  

Let $X$ be a projective variety over an algebraically closed field,
embedded by a (not necessarily complete) linear series. 
We begin by formulating a general version of the 
Minimal Resolution Conjecture (MRC), in analogy with the case of $\PP^n$,
predicting how the Betti numbers of 
a general subset of points of $X$ in the given embedding are related to the 
Betti numbers of $X$ itself. More precisely (cf. Theorem 
\ref{general_results} below),  
for a large enough general set of points $\Gamma$ on $X$, the Betti diagram 
consisting of the graded Betti numbers $b_{i,j}(\Gamma)$ is obtained from 
the Betti diagram of $X$ by adding two more nontrivial rows, at places 
well determined by the length of $\Gamma$. Recalling that the Betti diagram 
has the Betti number $b_{i,j}$ in the $(j,i)$-th position, and assuming 
that the two extra rows are indexed by $i=r-1$ and $i=r$, for some integer 
$r$, the MRC predicts that 
$$b_{i+1,r-1}(\Gamma)\cdot b_{i,r}(\Gamma)= 0,$$
i.e. at least one of the two Betti numbers on any "diagonal" is zero. 
As the difference $b_{i+1, r-1} - b_{i,r}$ can be computed exactly,
this implies a precise knowledge of the Betti numbers in these two rows.
Summing up, knowing the Betti diagram of $\Gamma$ would be the same as knowing 
the Betti diagram of $X$. A subtle question is however to understand
how the shape of the Betti diagram of $X$ influences whether MRC
is satisfied for points on $X$. An example illustrating this is given 
at the end of Section 1.

The Minimal Resolution Conjecture
 has been extensively studied in the case $X=\PP^n$. The
conjecture holds for $n\leq 4$ by results of
Gaeta, Ballico and Geramita, and Walter (see \cite{gaeta},
\cite{BG} and
\cite{walter}, respectively).
Moreover, Hirschowitz and Simpson 
proved in \cite{simpson}
 that it holds if the number of points is large enough with respect
to $n$. However, the conjecture does not hold in general:
it fails for every $n\geq 6$, $n\neq 9$ for almost $\sqrt{n}/2$
values of the number of points, by a result of Eisenbud, Popescu, Schreyer
and Walter (see \cite{epsw}). We refer to \cite{popescu} and
\cite{epsw} for a nice introduction and an account of the present status
of the problem in this case.

The main body of the paper is dedicated to a detailed study of MRC
in the case of curves. 
We will simply say that a curve \emph{satisfies MRC} in a given embedding 
if MRC is satisfied by a general set of points $\Gamma$ of any 
sufficiently large degree
(for the precise numerical statements see Section 1). We will also 
sometimes say that MRC holds for a line bundle $L$ if it 
holds for $C$ in the embedding given by $L$.  
Our main result says that MRC holds in the most significant case,
namely the case of canonically embedded nonhyperelliptic curves.

\begin{introtheorem}
If $C$ is a canonical curve, then $C$ satisfies MRC.
\end{introtheorem}

In contrast, under very mild assumptions on the genus, 
the MRC always fails in the case of curves of large degree, at well-determined 
spots in the Betti diagram (cf. Section 2 for precise details). The statement
$b_{2,r-1}\cdot b_{1,r}=0$, i.e. the case $i=1$, does hold though; this is 
precisely the Ideal Generation Conjecture, saying that the 
minimal number of generators of $I_\Gamma / I_X$ is as small as possible.

\begin{introtheorem}
(a) If $L$ is a very ample line bundle of degree $d\geq 2g$, then $IGC$ 
holds for $L$.
\newline
(b) If $g\geq 4$ and $L$ is a line bundle of degree $d\geq 2g+16$,
then there exists 
a value of $\gamma$ such that $C\subseteq\PP H^{0}(L)$ does not satisfy MRC 
for $i=\lfloor\frac{g+1}{2}\rfloor$. The same holds if $g\geq 15$ and
 $d\geq 2g+5$.
\end{introtheorem}

It is interesting to note that by the "periodicity" property of Betti
diagrams of 
general points on curves (see \cite{mustata} \S2), the theorem above
implies that 
on curves of high degree, MRC fails for sets of points of
arbitrarily large length.
This provides a very different picture from the case of projective space 
(cf. \cite{simpson}),
where asymptotically the situation is as nice as possible.

We explain the strategy involved in the proof of these results
 in some detail, as 
it appeals to some new geometric techniques in the study of
 syzygy related questions.
For simplicity we assume here that $C$ is a smooth curve embedded
 in projective space 
by means of a complete linear series corresponding to a very ample
 line bundle $L$
(but see \S2 for more general statements).
 A well-known geometric approach, developed  
by Green and Lazarsfeld in the study of syzygies of curves
 (see \cite{lazarsfeld} for a survey),
is to find vector bundle statements equivalent to the algebraic ones,
 via Koszul cohomology.
This program can be carried out completely in the case of MRC,
 and for curves we get a 
particularly clean statement. Assume that $M_L$ is the kernel of
 the evaluation map
$H^0(L)\otimes \OO_C\rightarrow L$ and $Q_L:=M_L^{*}$.
Then (cf. Corollary \ref{criterion} below) MRC holds 
for a collection of $\gamma\geq g$ general points on $C$
 if and only if the following 
is true:
$$h^0(\wedge^i M_L \otimes \xi)=0, ~{\rm for~all}~ i  ~{\rm and}~
\xi\in
 {\rm Pic}^{g-1+\lfloor{\frac{di}{n}}\rfloor}(C)~{\rm general}.\,\,\,\,(*)$$
Condition $(*)$ above is essentially the condition studied by Raynaud  
\cite{raynaud}, related to the existence of theta divisors for semistable 
vector bundles. In the particular situation of $\wedge^i M_L$, 
with $L$ a line bundle of large degree, it has been 
considered in \cite{popa} in order to produce base points for the determinant 
linear series on the moduli spaces $SU_C(r)$ of semistable bundles of rank $r$ 
and trivial determinant. A similar approach shows here the failure of
condition $(*)$ (and so of MRC) for $i=[\frac{g+1}{2}]$. On the
other hand, the fact that 
IGC holds is a
 rather elementary application of the Base Point Free Pencil Trick
\cite{ACGH} III \S3.

The case of canonical curves is substantially more involved, but in the end 
one is rewarded with a positive answer. As above, it turns out that MRC is 
equivalent to the vanishing:
$$h^0(\wedge^i Q\otimes \xi)=0, ~{\rm for~all}~i ~{\rm and~}~\xi\in
 {\rm Pic}^{g-2i-1}(C) 
~{\rm general},$$
where $Q$ is the dual of the bundle $M$ defined by the evaluation sequence: 
$$0\longrightarrow M\longrightarrow H^0(\omega_C)\otimes \OO_C\longrightarrow 
\omega_C\longrightarrow 0.$$
As the slope of $\wedge^i Q$ is $2i\in \ZZ$, this is in turn equivalent to the 
fact that $\wedge^i Q$ has
 a theta divisor $\Theta_{\wedge^i Q}\in {\rm Pic}^{g-2i-1}(C)$. 
On a fixed curve, if indeed a divisor,
$\Theta_{\wedge^i Q}$ will be identified as being precisely 
the difference variety $C_{g-i-1}-C_i \subseteq\Pici$ (cf. \cite{ACGH} Ch.V.D), where 
$C_n$ is the $n$-th symmetric product of $C$.
This is achieved via a 
filtration argument and a cohomology class calculation similar to the 
classical Poincar\'e theorem
 (cf. Proposition \ref{1}). A priori though, on an arbitrary curve 
the nonvanishing locus
 $\{\xi~|~h^0(\wedge^i Q\otimes \xi)\neq 0\}$ may be the whole $\Pici$, 
in which case this identification 
is meaningless. We overcome this problem by working with all curves at once, that 
is by setting up a similar universal construction on the moduli space 
of curves with marked points $\Mg1$.
 Here we slightly oversimplify the exposition in order 
to present the main idea,
 but for the precise technical details see Section 3. 
We essentially consider the "universal nonvanishing locus" in $\Mg1$: 
$$\cZ = \{(C,x_1,\ldots ,x_{g-i},y_1,\ldots ,y_{i+1})|h^0(\wedge^i
Q_C\otimes\OO(x_1+\ldots 
+x_{g-i}-y_1-\ldots -y_{i+1}))\neq 0\}.$$
The underlying idea is that the difference line bundles $\OO_C(x_1+\ldots 
+x_{g-i}-y_1-\ldots -y_{i+1})$ in fact cover the whole $\Pici$ (i.e. 
$C_{g-i}-C_{i+1}= \Pici$), and so 
for any given curve $C$,
 $\cZ|_C$ is precisely the nonvanishing locus described above. 
The advantage of writing
it in this form is that we are led to performing a computation on $\Mg1$ rather
than on a universal Picard,
 where for example one does not have a canonical choice of 
generators for the Picard group.  
A ``deformation to hyperelliptic
curves'' argument easily implies that MRC holds for general canonical 
curves, so $\cZ$ is certainly a divisor.
We then show that $\cZ$ is the degeneracy locus of a morphism of vector
bundles of the 
same rank and compute its class using a Grothendieck-Riemann-Roch 
argument (cf. Proposition \ref{degeneracy_class}). 

On the other hand, one can define an (a priori different) divisor $D$ in
$\Mg1$ which 
is a global analogue of the preimage of $C_{g-i-1}-C_i$ in 
$C^{g-i}\times C^{i+1}$ via
the difference map. It is convenient to see $D$ as the locus of curves 
with marked points $(C,x_1,\ldots ,x_{g-i},y_1,\ldots ,y_{i+1})$ having a
$\mathfrak g^1_g$ 
which contains $x_1,\ldots , x_{g-i}$ in a fiber and $y_1, \ldots ,
y_{i+1}$ in a different fiber.
An equivalent formulation of the discussion above is that 
$D\subseteq\cZ$, and in order for MRC to hold
 for all canonical curves one should 
have precisely $D = \cZ$. As we show that $D$
 is reduced (cf. Proposition \ref{reduced}), it suffices 
then to prove that the class of $D$ coincides with that of $\cZ$.
To this end we consider the closure of $D$ in the 
compactification $\overline{\cM}_{g,g+1}$,
 where the corresponding boundary condition is defined
by means of limit linear series. The computation of the class of $D$ via 
this closure is essentially independent of the rest of 
the paper. 
It relies on degeneration and enumerative techniques 
in the spirit of \cite{HM} and \cite{EH1}.

The results of both this and the computation
 of the class of $\mathcal{Z}$ are summarized 
in the following theorem. For the statement,
 we recall that ${\rm Pic}(\mathcal{M}_{g,n})_{\QQ}$
is generated by the class $\lambda$ of
 the Hodge bundle and the classes $\psi_j$, $1\leq j\leq n$,
where $\psi_j:= c_1(p_j^* \omega)$,
 with $\omega$ the relative dualizing sheaf on the universal 
curve $\mathcal{C}_g \rightarrow
 \mathcal{M}_g$ and $p_j : \Mgn \rightarrow \mathcal{C}_g$
the projection onto the $j$-th factor.

\begin{introtheorem}
The divisors $\cZ$ and $D$ defined above have the 
same class in $\Pic(\cM_{g,g+1})_{\QQ}$, namely
$$-\Bigl({g-1\choose i}-10{g-3\choose {i-1}}\Bigr)
\lambda+{{g-2}\choose i}\Psi_x
+{{g-2}\choose{i-1}}\Psi_y,$$
where $\Psi_x=\sum_{j=1}^{g-i}\psi_j$ and $\Psi_y=\sum_{j=g-i+1}^{g+1}\psi_j$. 
In particular $D = \cZ$.
\end{introtheorem}

As mentioned above, this implies that $\wedge^i Q$ always has a 
theta divisor, for all $i$, so equivalently that MRC holds
for an arbitrary canonical curve. We record the more 
precise identification of this theta divisor, which now follows
in general.

\begin{introcorollary}
For any nonhyperelliptic curve $C$, $\Theta_{\wedge^i Q} = C_{g-i-1}
- C_i$. 
\end{introcorollary}

In this particular form, our result answers positively a conjecture 
of R. Lazarsfeld.
It is worth mentioning that it also answers negatively a question that
was raised in connection with \cite{popa}, namely if $\wedge^i Q$ provide
base points for determinant linear series on appropriate moduli
spaces of vector bundles.

The paper is structured as follows. In Section 1 we give some 
equivalent formulations of the Minimal Resolution Conjecture and 
we describe the vector bundle setup used in the rest of the paper.
In Section 2 we treat the case of curves embedded with large degree,
proving IGC and showing that MRC fails. Section 3 is devoted to the 
main result, namely the proof of MRC for canonical curves, and here 
is where we look at the relationship with difference varieties and 
moduli spaces of curves with marked points. The divisor 
class computation in $\overline{\mathcal{M}}_{g,g+1}$, on which part of 
the proof relies, is carried out in Section 4 by means of limit linear
series.

\medskip
\noindent 
{\bf Acknowledgments.} We would like to thank D. Eisenbud and 
R. Lazarsfeld for sharing with us numerous useful ideas on the subject.

\section{\bf Several formulations of the Minimal Resolution Conjecture}

\subsection*{\bf Notations and conventions}

We work over an algebraically closed field $k$ which, unless
explicitly mentioned otherwise, has
arbitrary characteristic. Let $V$ be a vector space over $k$
with $\dim_kV=n+1$ and $S=Sym(V)\iso k[X_0,\ldots,X_n]$
the homogeneous coordinate ring of the corresponding
projective space $\PP V\iso\PP^n$. 

For a finitely generated graded $S$-module $N$,
the \emph{Betti~ numbers} $b_{i,j}(N)$ of $N$ are defined from the
minimal free resolution $F_{\bullet}$ of $N$ by
$$F_i=\oplus_{j\in\ZZ}S(-i-j)^{b_{i,j}(N)}.$$
The \emph{Betti ~diagram} of $N$ has in the $(j,i)$-th  position
the Betti number $b_{i,j}(N)$. The \emph{regularity}
${\rm reg}(N)$ of $N\neq 0$
can be defined as the index of the last nontrivial row
in the Betti diagram of $N$ (see \cite{eisenbud},
20.5 for
the connection with the cohomological definition).

We will use the computation of Betti numbers via
Koszul cohomology: $b_{i,j}(N)$ is the dimension over $k$
of the cohomology of the following piece of the
Koszul complex:
$$\wedge^{i+1}V\tensor N_{j-1}\longrightarrow
\wedge^iV\tensor N_j\longrightarrow\wedge^{i-1}V\tensor N_{j+1}$$
(see \cite{green} for details).

For an arbitrary subscheme $Z\subseteq \PP^n$, we denote by
$I_Z\subseteq S$ its saturated ideal and let $S_Z=S/I_Z$.
We denote by $P_Z$ and $H_Z$ the Hilbert polynomial and
Hilbert function of $Z$, respectively. The regularity 
${\rm reg}(Z)$ of $Z$ is defined to be the regularity of
$I_Z$, if $Z\neq \PP^n$, and $1$ otherwise. Notice that
with this convention, in the Betti diagram of $Z$,
which by definition is the Betti diagram of $S_Z$,
the last nontrivial row is always indexed by ${\rm reg}(Z)-1$. 

For a projective variety $X$, a line bundle $L$ on $X$, and
a linear series $V\subseteq H^0(L)$ which generates $L$, 
we denote by $M_V$ the vector bundle which is the kernel of the
evaluation map 
$$0\longrightarrow M_V \longrightarrow V\tensor\cO_X
 \overset{ev}{\longrightarrow} 
L\longrightarrow 0.$$
When $V=H^0(L)$ we use the
notation $M_L:=M_V$. If $C$ is a smooth curve of genus $g\geq 1$,
and $\omega_C$ is the canonical line bundle, then 
$M_C$ denotes the vector bundle $M_{\omega_C}$. The dual vector bundles
will be denoted by $Q_V$, $Q_L$ and $Q_C$, respectively. 
Whenever there is no risk of confusion, we will simply write
$M$ and $Q$, instead of $M_C$ and $Q_C$.

\subsection*{\bf The Minimal Resolution Conjecture for points on
 embedded varieties.} 
In this section $X\subseteq\PP\,V\iso\PP^n$ is a fixed irreducible projective 
variety of positive dimension.
 We study the Betti numbers of a general set of $\gamma$ points $\Gamma
\subseteq X$.
 Since the Betti numbers are upper semicontinuous functions,
for every positive integer $\gamma$, there is an open subset $U_{\gamma}$ of
$X^{\gamma}\setminus\cup_{p\neq q}\{x\colon x_p=x_q\}$ such that
for all $i$ and $j$, $b_{i,j}(\Gamma)$ takes its minimum value
for $\Gamma\in U_{\gamma}$. Notice that as the regularity is bounded in terms
of $\gamma$, we are concerned with finitely many Betti numbers.
From now on, $\Gamma$ \emph{general} means $\Gamma\in U_{\gamma}$.

It is easy to determine the Hilbert function of a general set of points 
$\Gamma$ 
in terms of the Hilbert function of $X$ 
(see \cite{mustata}). We have the following:

\begin{proposition}\label{hilbert_function}
If $\Gamma\subseteq X$ is a general set of $\gamma$ points, then
$$H_{\Gamma}(t)={\rm min}\,\{H_X(t),\gamma\}.$$
\end{proposition}

To determine the Betti numbers of a general set of points $\Gamma$ is a 
much more subtle problem. If $\gamma$ is large enough, then 
the Betti diagram of $\Gamma$ looks as follows: in the upper part we
have the Betti diagram of $X$ and there are two extra nontrivial rows at
the bottom. Moreover, the formula in Proposition~\ref{hilbert_function}
gives an expression for the differences of the Betti numbers in these
last two rows. We record the formal statement in the following theorem
and for the proof we refer to \cite{mustata}.

\begin{theorem}\label{general_results}
Assume that $\Gamma\subseteq X$ is a general set of $\gamma$ points,
with $P_X(r-1)\leq\gamma<P_X(r)$ for some $r\geq m+1$, where $m={\rm reg}\,X$.

\item{\rm (i)} For every $i$ and $j\leq r-2$, we have $b_{i,j}(\Gamma)
=b_{i,j}(X)$.
\item{\rm (ii)} $b_{i,j}(\Gamma)=0$, for $j\geq r+1$ and there is an  
$i$ such that $b_{i,r-1}(\Gamma)\neq 0$.
\item{\rm (iii)} For every $j\geq m$, we have
$$b_{i,j}(\Gamma)=b_{i-1,j+1}(I_{\Gamma}/I_X)=b_{i-1,j+1}(\oplus_{l\geq 0}
H^0(\cI_{\Gamma/X}(l))).$$
\item{\rm (iv)} If $d=\dim\,X$, then for every $i\geq 0$, we have
$b_{i+1,r-1}(\Gamma)-b_{i,r}(\Gamma)=Q_{i,r}(\gamma)$, where
$$Q_{i,r}(\gamma)=\sum_{l=0}^{d-1}
(-1)^l{{n-l-1}\choose{i-l}}\Delta^{l+1}P_X(r+l)-
{n\choose i}(\gamma-P_X(r-1)).$$
\end{theorem}

We will focus our attention on the Betti numbers in the bottom two rows
in the Betti diagram of $\Gamma$. The equation in Theorem~\ref
{general_results} (iv) gives lower bounds for these numbers, namely
$b_{i+1,r-1}(\Gamma)\geq {\rm max}\,\{Q_{i,r}(\gamma),0\}$
and $b_{i,r}(\Gamma)\geq {\rm max}\,\{-Q_{i,r}(\gamma),0\}$.

\medskip
\begin{definition}
In analogy with the case $X=\PP^n$
(see \cite{lorenzini1} and \cite{lorenzini2}), we say that the 
\emph{Minimal Resolution Conjecture} (to which we refer from now on as MRC)
holds for a fixed value of $\gamma$ as above
if for every $i$ and every  general set $\Gamma$,
 $b_{i+1,r-1}(\Gamma)={\rm max}\,\{Q_{i,r}(\gamma),0\}$
 and $b_{i,r}(\Gamma)={\rm max}\,\{-Q_{i,r}(\gamma),0\}$.
 Equivalently, it says that 
$$b_{i+1,r-1}(\Gamma)\cdot b_{i,r}(\Gamma)
=0 {\rm ~for~ all~} i.$$
\end{definition}

This conjecture has been extensively studied in the case $X=\PP^n$,
 $L=\OO_{\PP^n}(1)$. 
It is known to hold for small values of $n$ ($n=2$, $3$ or $4$) and for
large values of $\gamma$, depending on $n$, but not in
general. In fact, it has been shown that for every $n\geq 6$, $n\neq 9$,
MRC fails for almost $\sqrt{n}/2$ values of $\gamma$ (see
\cite{epsw}, where one can find also a detailed account of the problem).

Note that the assertion in MRC holds obviously for $i=0$. The first
nontrivial case $i=1$ is equivalent by Theorem~\ref{general_results}
to saying that the minimal number of generators of $I_{\Gamma}/I_X$
is as small as possible. This suggests the following:

\medskip
\begin{definition}
We say that the \emph{Ideal Generation Conjecture}
(IGC, for short) holds for
$\gamma$ as above if for a general set of points $\Gamma\subseteq X$ of
 cardinality $\gamma$,
we have $b_{2,r-1}(\Gamma)\cdot b_{1,r}(\Gamma)=0$.
\end{definition}

\smallskip

\begin{example}(\cite{mustata})\label{case1}
 MRC holds for every $X$ when $\gamma=P_X(r-1)$,
since in this case $b_{i,r}(\Gamma)=0$ for every $i$.
Similarly, MRC holds for every $X$ when $\gamma= P_X(r)-1$,
since in this case $b_{1,r-1}=1$ and $b_{i,r-1}(\gamma)=0$ for $i\geq 2$.
\end{example}

We derive now a cohomological interpretation of MRC.
From now on we assume that $X$ is nondegenerate, so that
we have $V\subseteq H^0(\cO_X(1))$.
Using a standard Koszul cohomology argument, we can express the
Betti numbers in the last two rows of the Betti diagram of $\Gamma$
as follows.

\begin{proposition}\label{cohom}
With the above notation, we have in general for every $i\geq 0$ 
$$b_{i+1,r-1}(\Gamma)=h^0(\wedge^i M_V\tensor\cI_{\Gamma/X}(r)),$$
$$b_{i,r}(\Gamma)=h^1(\wedge^i M_V\tensor\cI_{\Gamma/X}(r)).$$
\end{proposition}

\begin{proof}
We compute the Betti numbers via Koszul cohomology,
using the formula in Theorem~\ref{general_results} (iii).

Consider the complex:
$$ \wedge^iV\tensor H^0(\I_{\Gamma/X}(r))\overset{f}{\longrightarrow}
\wedge^{i-1}V\tensor H^0(\I_{\Gamma/X}(r+1))
\overset{h}{\longrightarrow}\wedge^{i-2}V\tensor H^0(\I_{\Gamma/X}(r+2))$$

Since $H^0(\cI_{\Gamma/X}(r-1))=0$, it follows
that $\dim_k({\rm Ker}f)
=b_{i+1,r-1}(\Gamma)$ and $\dim_k({\rm Ker}\,h/{\rm Im}\,f)
=b_{i,r}(\Gamma)$.
 The exact sequence
$$0\longrightarrow M_V\longrightarrow V\tensor\OO_X\longrightarrow\OO_X(1)
\longrightarrow 0$$ induces long exact sequences
$$0\longrightarrow\wedge^i M_V\longrightarrow\wedge^i V\tensor\OO_X
\longrightarrow\wedge^{i-1} M_V\tensor\OO_X(1)\longrightarrow 0.\,\,\,\,
(*)$$
By tensoring with $\cI_{\Gamma/X}(r)$ and taking global sections,
we get the exact sequence
$$H^0(\wedge^i M_V\tensor\I_{\Gamma/X}(r))
\hookrightarrow\wedge^iV\tensor H^0(\cI_{\Gamma/X}(r))
\overset{f}{\longrightarrow}\wedge^{i-1}V\tensor H^0(\cI_{\Gamma/X}(r+1)).$$
This proves the first assertion in the proposition.

\smallskip

We have a similar exact sequence:
$$H^0(\wedge^{i-1} M_V\tensor\cI_{\Gamma/X}(r+1))
\hookrightarrow\wedge^{i-1} V\tensor H^0(\cI_{\Gamma/X}(r+1))
\overset{h}{\longrightarrow}
\wedge^{i-2}V\tensor H^0(\cI_{\Gamma/X}(r+2)).$$
Therefore $b_{i,r}(\Gamma)$ is the dimension
over $k$ of the cokernel of
$$g\,:\,\wedge^iV\tensor H^0(\cI_{\Gamma/X}(r))\longrightarrow
H^0(\wedge^{i-1}M_V\otimes\cI_{\Gamma/X}(r+1)).$$

Using again the exact sequence $(*)$, by tensoring with $\cI_{\Gamma/X}(r)$
and taking a suitable part of the long exact sequence, we get:
$$\wedge^iV\tensor H^0(\cI_{\Gamma/X}(r))\longrightarrow
H^0(\wedge^{i-1} M_V\tensor\cI_{\Gamma/X}(r+1))\longrightarrow$$
$$H^1(\wedge^i M_V\tensor\cI_{\Gamma/X}(r))\longrightarrow\wedge^iV
\tensor H^1(\cI_{\Gamma/X}(r)).$$

Since ${\rm reg}\,\Gamma\leq r+1$, we have ${\rm reg}\,\cI_{\Gamma/X}\leq r+1$
and therefore
$H^1(\cI_{\Gamma/X}(r))=0$. From the above exact sequence we see
that ${\rm Coker}\,g\simeq
H^1(\wedge^i M_V\tensor\cI_{\Gamma/X}(r))$, which proves
the second assertion of the proposition.
\end{proof}

\smallskip

\begin{remark}\label{higher}
 The higher cohomology groups 
$H^p(\wedge^i M_V\tensor\I_{\Gamma/X}(r))$, $p\geq 2$, always vanish. 
Indeed, using the exact sequences in the proof of the proposition,
we get $$h^p(\wedge^i M_V\tensor\I_{\Gamma/X}(r))=
h^1(\wedge^{i-p+1} M_V\tensor\I_{\Gamma/X}(r+p-1))=
b_{i-p+1,r+p-1}(\Gamma)=0.$$
Therefore we have $Q_{i,r}(\gamma)=\chi(\wedge^i M_V\tensor\I_{\Gamma/X}(r))$
and MRC can be interpreted as saying that for general $\Gamma$,
the cohomology of $\wedge^i M_V\tensor\I_{\Gamma/X}(r)$ is supported
in cohomological degree either zero or one.
\end{remark}

\smallskip

In the case of a curve $C$, 
MRC can be reformulated using 
Proposition~\ref{cohom}
in terms of general line bundles on $C$.
We will denote by $\lfloor x\rfloor$ and $\lceil x\rceil$
the integers defined by $\lfloor x\rfloor\leq x <\lfloor x\rfloor +1$
and $\lceil x \rceil -1<x\leq\lceil x\rceil$. 

\begin{corollary}\label{criterion}
Suppose that $C\subseteq\PP\,V$
is a nondegenerate, integral curve of arithmetic genus  $g$ 
and degree $d$.
We consider the following two statements:

\item{\rm (i)} For every  $i$ and for a general line bundle
$\xi\in \Pic^j(C)$, where $j=g-1+\lceil{\frac{di}{n}}\rceil$, we have
$H^1(\wedge^i M_V\tensor \xi)=0$. 
\item{\rm (ii)} For every $i$ and for a general line bundle
$\xi\in\Pic^j(C)$, where $j=g-1+\lfloor{\frac{di}{n}}\rfloor$, we have
$H^0(\wedge^i M_V\tensor \xi)=0$.

Then MRC holds for $C$ for every $\gamma\geq \max\{g, P_C({\rm reg}\,X)\}$
 if and only if
both {\rm (i)} and {\rm (ii)} are true. Moreover, if $C$ is locally Gorenstein,
then {\rm (i)} and {\rm (ii)} are equivalent.
\end{corollary}

\begin{proof} 
If $\gamma\geq g$, then for a general set $\Gamma$ of $\gamma$ points,
$\I_{\Gamma/C}$ is a general line bundle on $C$ of degree $-\gamma$.
Since in this case $\I_{\Gamma/C}(r)$ is a general line bundle of
degree $j=dr-\gamma$ and $d(r-1)+1-g\leq \gamma\leq dr+1-g$,
Proposition~\ref{cohom}
says that MRC holds for every $\gamma\geq \max\{g, P_C({\rm reg}\,C)\}$
 if and only if
for every $j$ such that $g-1\leq j\leq d+g-1$ and for
a general line bundle $\xi'\in\Pic^j(C)$,
either $H^0(\wedge^i M_V\tensor \xi')=0$ or
$H^1(\wedge^i M_V\tensor \xi')=0$. 

Since $\dim\,C=1$, we have $b_{i+1,r-1}(\Gamma)-b_{i,r}(\Gamma)=
d{n-1\choose i} -(\gamma-P_C(r-1)){n\choose i}$. It follows immediately
that $b_{i+1,r-1}(\Gamma)-b_{i,r}(\Gamma)\geq 0$ if and only if
$j\geq g-1+di/n$.

The first statement of the corollary follows now from the fact that
if $E$ is a vector bundle on a curve and $P$ is a point, then
$H^0(E)=0$ implies $H^0(E\tensor\OO(-P))=0$ and
$H^1(E)=0$ implies $H^1(E\tensor\OO(P))=0$. The last statement
follows from Serre duality and the isomorphism 
$\wedge^iQ_V\simeq\wedge^{n-i}M_V\otimes\cO_C(1)$.
\end{proof}

\begin{remark}\label{caseIGC}
 The corresponding assertion for IGC says that
$X$ satisfies IGC for every $\gamma\geq \max\{g, P_C({\rm reg}\,C)\}$
 if and only if both 
(i) and (ii) are true for $i=1$. Note that if $X$ is locally Gorenstein,
then by Serre duality condition (ii) for $i=1$ is equivalent to
condition (i) for $i=n-1$.
\end{remark}

\begin{remark}\label{half}
If $C$ is a locally  Gorenstein integral curve
such that $d/n\in\ZZ$, then in order to check MRC for all
$\gamma\geq \max\{g, P_C({\rm reg}\,C)\}$,
 it is enough to check condition (i) in Corollary~\ref
{criterion} only for $i\leq n/2$. Indeed, using Serre duality
and Riemann-Roch, we see that the conditions for $i$
and $n-i$ are equivalent.
\end{remark}

In light of Corollary~\ref{criterion}, we make the following:

\begin{definition}
If $C\subseteq \PP\,V$ is a nondegenerate integral curve of arithmetic genus
$g$ and regularity $m$, we say that $C$ satisfies MRC if 
a general set of $\gamma$ points on $C$ satisfies MRC
for every $\gamma\geq{\rm max}\{g, P_C(m)\}$. If $L$ is a very ample line
bundle on a curve $C$ as before, we say that $L$ satisfies MRC if
$C\subset\PP\,H^0(L)$ satisfies MRC. Analogous definitions are made for IGC.
\end{definition}

\begin{example}({\bf Rational quintics in $\PP^3$.})\label{quintics}
We illustrate the above discussion in the case of smooth rational quintic
curves in $\PP^3$. We consider two explicit examples, the first when the curve
lies on a (smooth) quadric and the second when it does not.
Let $X$ be given parametrically by
$(u,v)\in\PP^1\longrightarrow (u^5,u^4v,uv^4,v^5)\in\PP^3$, so that
it lies on the quadric $X_0X_3=X_1X_2$. The Betti diagram of $X$ is
\begin{center}
\renewcommand{\arraystretch}{1.25}
\begin{tabular}[l]{c|c}
\hline
0 & 1 -- -- -- \\
1 & -- 1 -- --  \\
2 & -- -- -- -- \\
3 & -- 4 6 2  \\
\end{tabular}
\end{center}
and if $\Gamma\subset X$ is a set of $28$ points, then the Betti diagram of
$\Gamma$ is
\begin{center}
\renewcommand{\arraystretch}{1.25}
\begin{tabular}[l]{c|c}
\hline
0 & 1 -- -- -- \\
1 & -- 1 -- -- \\
2 & -- -- -- -- \\
3 & -- 4 6 2 \\
4 & -- -- -- -- \\
5 & -- 3 4 1 \\
6 & -- -- 2 2 \\
\end{tabular}
\end{center}
As $b_{3,5}(\Gamma)=1$ and $b_{2,6}(\Gamma)=2$, we see that MRC
is not satisfied by $X$ for this number of points.

Let now $Y$  be the curve given parametrically by
$(u,v)\in\PP^1\longrightarrow (u^5+u^3v^2, u^4v-u^2v^3,uv^4,v^5)\in\PP^3$.
In this case $Y$ does not lie on a quadric, and in fact, its
Betti diagram is given by
\begin{center}
\renewcommand{\arraystretch}{1.25}
\begin{tabular}[l]{c|c}
\hline
0 & 1 -- -- -- \\
1 & -- -- -- -- \\
2 & -- 4 3 -- \\
3 & -- 1 2 1 \\
\end{tabular}
\end{center}
If $\Gamma'\subset Y$ is a set of $28$ points, then the Betti diagram of
$\Gamma'$ is
\begin{center}
\renewcommand{\arraystretch}{1.25}
\begin{tabular}[l]{c|c}
\hline
0 & 1 -- -- -- \\
1 & -- -- -- -- \\
2 & -- 4 3 -- \\
3 & -- 1 2 1 \\
4 & -- -- -- -- \\
5 & -- 3 4 -- \\
6 & -- -- 1 2 \\
\end{tabular}
\end{center}
which shows that MRC is satisfied for $Y$ and this number of points.

These two examples show the possible behavior with respect to the
MRC for smooth rational quintics in $\PP^3$. The geometric condition
of lying on a quadric translates into a condition on the splitting type
of $M_V=\Omega_{\PP^3}(1)\vert_X$. More precisely, it is proved in
\cite{ev} that if $X\subset\PP^3$ is a smooth rational quintic curve, then
$X$ lies on a quadric if and only if we have
$\Omega_{\PP^3}(1)\vert_X\simeq\cO_{\PP^1}(-3)\oplus\cO_{\PP^1}(-1)^{\oplus 2}$
(the other possibility, which is satisfied by
a general such quintic, is that $\Omega_{\PP^3}(1)\vert_X
\simeq\cO_{\PP^1}(-1)\oplus\cO_{\PP^1}(-2)^{\oplus 2}$).
Corollary~\ref{criterion} explains therefore the behaviour with respect
to MRC in the above examples.
\end{example}

\section{\bf Curves of large degree and a counterexample to MRC}

In this section we assume that $C$ is a smooth projective curve
of genus $g$ and $L$ is a very ample line bundle on $C$. Our aim is to 
investigate whether $C$ satisfies MRC, or at least IGC,  for every
$\gamma\geq g$, in the embedding given by the complete linear series $|L|$. 
As before, $m$ will denote the regularity of $C$.

\begin{example}\label{g_small}
If $g=0$ or $1$, then $C$ satisfies MRC for all $\gamma\geq P_C(m)$ 
in every embedding given by a complete linear series
(see \cite{mustata}, Proposition 3.1).
\end{example}

In higher genus we will concentrate on the study
 of MRC for canonical curves and 
curves embedded with high degree, in direct analogy with the syzygy questions 
of Green-Lazarsfeld (cf. \cite{gl1}, \cite{gl2}, \cite{gl3}).
The main conclusion of this section will be that, while 
IGC is satisfied in both situations, the high-degree embeddings
always fail to satisfy MRC at a well-specified spot in the Betti diagram. 
This is in contrast with our main result, proved in \S3, that MRC 
always holds for canonical curves, and the arguments involved here  
provide an introduction to that section. The common theme of the proofs is the
vector bundle interpretation of MRC described in \S1.

\medskip
\noindent 
{\bf Review of filtrations for $Q_L$ and $Q$ \cite{lazarsfeld}.}
 Here we recall a basic property of 
the vector bundles $Q_L$ which will be essential for our arguments. 
Let $L$ be a very ample line bundle on $C$
of degree $d$, and recall from \S1 that $Q_L$ is given by the defining sequence
$$0\longrightarrow L^{-1}\longrightarrow H^0(L)^{*}\otimes
 \OO_{C}\longrightarrow
Q_{L}\longrightarrow 0.$$
Assume first that $L$ is non-special and $x_{1},\ldots,x_{d}$ are the points
of a general hyperplane section of $C\subseteq\PP H^{0}(L)$.
 One shows (see e.g. \cite{lazarsfeld} \S1.4) 
that there exists an exact sequence 
\begin{equation}\label{large_degree_filtration}
0\longrightarrow \underset{i\in \{1,\ldots, d-g-1\}}
{\bigoplus}\OO_{C}(x_{i})\longrightarrow Q_{L}
\longrightarrow \OO_{C}(x_{d-g}+ \ldots +x_{d})\longrightarrow 0.
\end{equation}
On the other hand, assuming that $C$ is nonhyperelliptic and $L=\omega_C$, if 
$x_{1},\ldots,x_{2g-2}$ are the points of a general hyperplane section,
the analogous sequence reads:
\begin{equation}\label{canonical_filtration}
0\longrightarrow \underset{i\in \{1,\ldots, g-2\}}
{\bigoplus}\OO_{C}(x_{i})\longrightarrow Q
\longrightarrow \OO_{C}(x_{g-1}+ \ldots +x_{2g-2})\longrightarrow 0.
\end{equation}

\medskip
We start by looking at the case of curves embedded with large degree.
 The main results 
are summarized in the following:

\begin{theorem}\label{curves_of_large_degree}
{\rm (a)} If $L$ is a very ample line bundle of degree $d\geq 2g$, then $IGC$ 
holds for $L$.
\newline
{\rm (b)} If $g\geq 4$ and $L$ is a line bundle of degree $d\geq 2g+10$,
 then there exists 
a value of $\gamma$ such that $C\subseteq\PP H^{0}(L)$ does not satisfy MRC 
for $i=\lfloor\frac{g+1}{2}\rfloor$. The same holds if $g\geq 14$ and
 $d\geq 2g+5$.
\end{theorem}

\begin{proof}
(a)
Let $L$ be a very ample line bundle of degree $d\geq 2g$.
By Corollary \ref{criterion} and Serre duality, it is easy to see that IGC 
holds for $L$ if:
\medskip
\newline
(i) $h^{1}(Q_{L}\otimes \eta)=0$ for $\eta\in {\rm Pic}^{g-2}(C)$ general
\newline
and
\newline
(ii) $h^{0}(Q_{L}\otimes \eta)=0$ for $\eta\in {\rm Pic}^{g-3}(C)$ general.

Condition (i) is a simple consequence of the filtration
 (\ref{large_degree_filtration}). 
More precisely, if $x_{1},\ldots,x_{d}$ are the points
of a general hyperplane section of $C\subseteq\PP H^{0}(L)$,
 from the exact sequence
$$0\longrightarrow \underset{i\in \{1,\ldots, d-g-1\}}
{\bigoplus}\OO_{C}(x_{i})\longrightarrow Q_{L}
\longrightarrow \OO_{C}(x_{d-g}+ \ldots +x_{d})\longrightarrow 0.$$
we conclude that it would be enough to prove:
$$h^{1}(\eta(x_{i}))=0 {\rm ~and~} h^{1}(\eta(x_{d-g}+ \ldots +x_{d}))=0$$
for $\eta\in {\rm Pic}^{g-2}(C)$ general. Now for every
 $i\in \{1,\ldots ,d-g-1\}$, $\eta(x_{i})$
is a general line bundle of degree $g-1$, so $h^{1}(\eta(x_{i}))=0$.
 On the other hand
${\rm deg}~\eta(x_{d-g}+ \ldots +x_{d})\geq2g-1$, so clearly
 $h^{1}(\eta(x_{d-g}+ \ldots +x_{d}))=0$.

For condition (ii) one needs a different argument.
 By twisting the defining sequence of $Q_{L}$: 
$$0\longrightarrow L^{-1}\longrightarrow
H^{0}(L)^{*}\otimes \OO_{C}\longrightarrow
Q_{L}\longrightarrow 0$$
by $\eta\in {\rm Pic}^{g-3}(C)$ general and taking cohomology,
 we see that (ii) holds if and only 
if the map
$$\alpha^{*}:H^{1}(L^{-1}\otimes \eta)\rightarrow
H^{0}(L)^{*}\otimes H^{1}(\eta)$$
is injective, or dually if and only if the cup-product map
$$\alpha: H^{0}(L)\otimes H^{0}(\omega_{C}
\otimes \eta^{-1})\rightarrow H^{0}(L\otimes \omega_{C}
\otimes \eta^{-1})$$
is surjective. We make the following:
\medskip
\newline
\emph{Claim.} $|\omega_{C}\otimes \eta^{-1}|$
 \emph{is a base point free pencil}.

\medskip
Assuming this for the time being, one can apply the
 Base Point Free Pencil Trick
(see \cite{ACGH} III \S3) to conclude that
$${\rm Ker~}\alpha= H^{0}(L\otimes \omega_{C}^{-1}\otimes \eta).$$
But $L\otimes \omega_{C}^{-1}\otimes \eta$ is a general line bundle of
 degree $d-g-1\geq g-1$  
and so $h^{1}(L\otimes \omega_{C}^{-1}\otimes \eta)=0$.
 By Riemann-Roch this means
$h^{0}(L\otimes \omega_{C}^{-1}\otimes \eta)=d-2g$.
 On the other hand $h^{0}(L)=d-g+1$,
$h^{0}(\omega_{C}\otimes \eta^{-1})=2$ and
 $h^{0}(L\otimes \omega_{C}\otimes \eta^{-1})=d+2$, so
$\alpha$ must be surjective.

We are only left with proving the claim. Since $\eta\in {\rm Pic}^{g-3}(C)$
 is general,
$h^{0}(\eta)=0$, and so we easily get:
$$ h^{0}(\omega_{C}\otimes \eta^{-1})=h^{1}(\eta)=g-1-(g-3)=2.$$
Also, for every $p\in C$, $\eta(p)\in {\rm Pic}^{g-2}(C)$ is general,
 hence still noneffective. Thus:
$$h^{0}(\omega_{C}\otimes \eta^{-1}(-p))=h^{1}(\eta(p))=g-1-(g-2)=1.$$
This implies that $|\omega_{C}\otimes \eta^{-1}|$ is base point free.

\medskip
(b) Here we follow an argument in \cite{popa} leading to the required
nonvanishing statement.
First note that it is clear from (\ref{large_degree_filtration})
 that for every $i$ with $1\leq i\leq d-g-1$
there is an inclusion
$$\OO_{C}(x_{1}+\ldots +x_{i})\hookrightarrow \wedge^{i}Q_{L},$$
where $x_{1},\ldots , x_{i}$ are general points on $C$.
This immediately implies that 
$$h^{0}(\wedge^{i}Q_{L}\otimes \OO_{C}(E_{i}-D_{i}))\neq 0,$$
where $E_{i}$ and $D_{i}$ are general effective divisors on $C$ of degree $i$.
On the other hand 
we use the fact (see e.g. \cite{ACGH} Ex. V. D) 
that every line bundle $\xi \in {\rm Pic}^{0}(C)$ can be written as a 
difference
$$\xi=\OO_{C}(E_{\lfloor\frac{g+1}{2}\rfloor}-
D_{\lfloor\frac{g+1}{2}\rfloor}),$$
which means that 
$$h^{0}(\wedge^{\lfloor\frac{g+1}{2}\rfloor}
Q_{L}\otimes \xi)\neq 0,~\forall \xi 
\in {\rm Pic}^{0}(C) {\rm ~general~}.$$
Now by Serre duality:
$$H^{0}(\wedge^{i}Q_{L}\otimes \xi)\cong H^{1}(\wedge^{i}M_{L}\otimes 
\omega_{C}\otimes \xi^{-1})^{*},$$
so that Corollary \ref{criterion} easily implies that $C$
 does not satisfy MRC for 
$i=\lfloor\frac{g+1}{2}\rfloor$ as long as $2g-2\geq g-1+ \frac{di}{d-g}$.
A simple computation gives then the stated conclusion.
\end{proof}

\begin{remark}
Motivation for the argument in (b) above was quite surprisingly
provided by the study \cite{popa} of 
the base locus of the determinant linear series on the moduli space $SU_C(r)$ 
of semistable bundles of rank $r$ and trivial determinant on a curve $C$. 
In fact this argument produces 
explicit base points for the determinant linear series under appropriate
 numerical conditions.
\end{remark}

\begin{remark}
The technique in Theorem \ref{curves_of_large_degree} (b) can be extended
 to produce 
examples of higher dimensional varieties for which appropriate
 choices of $\gamma$ 
force the failure of MRC for general sets of $\gamma$ points.
 More precisely, the 
varieties in question are projective bundles $\PP E \rightarrow C$ over a
 curve $C$, 
associated to very ample 
vector bundles $E$ on $C$ of arbitrary rank and large degree,
 containing sub-line bundles of 
large degree. Using the interpretation given in Proposition \ref{cohom},
 the problem 
is reduced to a cohomological question about the exterior powers
 $\wedge^i M_E$, where 
$M_E$ is defined analogously as the kernel of the evaluation map
$$0\longrightarrow M_E \longrightarrow H^0(E)\otimes
 \OO_X\longrightarrow E\longrightarrow 0.$$
This question is then treated essentially as above,
 and we do not enter into details. 
Unfortunately once a bundle $E$ of higher rank is fixed,
 this technique does not seem 
to produce couterexamples for arbitrarily large values of $\gamma$,
 as in the case of 
line bundles. Such examples would be very interesting,
 in light of the asymptotically nice
behavior of general points in $\PP^n$ (cf. \cite{simpson}).
\end{remark}

Finally we turn to the case of canonical curves with the goal of
 providing an introduction 
to the main result in Section 3. Let $C$ be 
a nonhyperelliptic curve of genus $g$, $V=H^0(\omega_C)$
and $C\hookrightarrow\PP V\iso\PP^{g-1}$ the canonical embedding.
We note here that an argument similar to
 Theorem \ref{curves_of_large_degree} (a) 
immediately implies IGC for $C$. This will be later subsumed in the general 
Theorem \ref{canonical_curves}.  

\begin{proposition}\label{igc_canonical}
IGC holds for the canonical curve $C$.
\end{proposition}
\begin{proof}
The argument is similar (and in fact simpler) to the proof of (ii) in 
Theorem \ref{curves_of_large_degree} (a). 
In this case, again by interpreting Proposition~\ref{cohom}
(se Remark~\ref{caseIGC}
 IGC holds if and only if
$$H^{0}(Q\otimes \xi)=0~{\rm  for}~\xi\in {\rm Pic}^{g-3}(C) ~{\rm general}.$$
This is in turn equivalent to the surjectivity of the multiplication map:
$$H^{0}(\omega_{C})\otimes H^{0}(\omega_{C}\otimes \xi^{-1})\rightarrow
H^{0}(\omega_{C}^{\otimes 2}\otimes \xi^{-1}),$$
which is again a quick application  of the Base Point Free Pencil Trick.
\end{proof}

The geometric picture in the present case of canonical curves
 can be described a
little more precisely. In fact, for $\xi\in {\rm Pic}^{g-3}(C)$, we have
$$\mu(Q\otimes \xi)=g-1,$$
where  $\mu(E):={\rm deg}\,(E)/{\rm rk}\,(E)$
denotes in general the \emph{slope} of the vector bundle $E$.
By standard determinantal results, the subset 
$$\Theta_{Q}:=\{\xi~|~ h^{0}(Q\otimes \xi)\neq 0\}\subseteq{\rm Pic}^{g-3}(C)$$
is either a divisor or the whole variety.
The statement of IGC is then equivalent to saying that  
$\Theta_{Q}$ is indeed a divisor
in ${\rm Pic}^{g-3}(C)$ (one says that $Q$ \emph{has a theta divisor}). 
A simple filtration argument based on the sequence (3) above shows that in fact
$$\Theta_{Q}=C_{g-2}-C:=\{\OO_{C}(p_{1}+\ldots+p_{g-2}-q)~|~p_{1},\ldots ,
p_{g-2}, q \in C\},$$
which has already been observed by Paranjape and Ramanan in
\cite{paranjape}
A generalization of this observation to the higher exterior powers
 $\wedge^iQ$ will be 
the starting point for our approach to proving MRC for canonical curves in
 what follows.

\section{\bf MRC for canonical curves}

In this section $C$ will be a canonical curve, i.e. a smooth curve
of genus $g$ embedded in $\PP^{g-1}$ by the canonical linear series
$|\omega_C|$ (in particular $C$ is not hyperelliptic). Our goal is to prove
the following:

\begin{theorem}\label{canonical_curves}
If $C$ is a canonical curve, then $C$ satisfies MRC.
\end{theorem}

\begin{remark}
In fact, since $C$ is canonically embedded, its regularity is
$m=4$, and as $g\geq 3$ we always have $P_C(m)=7(g-1)\geq g$.
Thus the statement means that MRC holds for every $\gamma \geq P_C(m)$.
\end{remark}

The general condition required for a curve to 
satisfy MRC which was stated in Corollary \ref{criterion}
(see also Remark~\ref{half}) 
takes a particularly clean form in the case of canonical embeddings.
We restate it for further use.

\begin{lemma}\label{restatement}
Let $C$ be a canonical curve. Then $C$ satisfies MRC if and only if, for all 
$1 \leq i\leq \frac{g-1}{2}$ we have
$$h^0(\wedge^i M\otimes \eta) =
h^1(\wedge^i M\otimes \eta)=0, ~{\rm for}~ \eta\in 
{\rm Pic}^{g+2i-1}(C)~{\rm general},$$ 
or equivalently
$$(*) ~~h^0(\wedge^i Q\otimes \xi) =h^1(\wedge^i Q\otimes \xi)=0, 
~{\rm for}~ \xi\in \Pici~
{\rm general}.$$ 
\end{lemma} 

\begin{remark}
Note that $\mu(Q) = 2$, so $\mu(\wedge^i Q)=2i\in \ZZ$.
 This means that the condition 
$(*)$ in Lemma~\ref{restatement} is equivalent to saying that $\wedge^i Q$
 has a theta divisor (in 
$\Pici$), which we denote $\Theta_{\wedge^i Q}$.
 In other words, the set defined by 
$$\Theta_{\wedge^i Q}:= \{ \xi\in \Pici~|~h^0(\wedge^i Q\otimes \xi)\neq 0\}$$
with the scheme structure
 of a degeneracy locus of a map of vector bundles of the same 
rank is an actual divisor as expected (cf. \cite{ACGH} II \S4).
\end{remark}

\subsection*{\bf Hyperelliptic curves.}  Note that the statement $(*)$
in Lemma~\ref{restatement} makes sense even for hyperelliptic curves.
 Again $Q$ is the dual of $M$, where $M$ is the 
kernel of the evaluation map for the canonical line
bundle.
Therefore we will say slightly abusively 
that MRC is satisfied for some smooth curve of genus
$g\geq 2$ if $(*)$ is satisfied for all $i$, $1\leq i\leq (g-1)/2$.
In fact, the hyperelliptic case is the only one for which we can give a direct
argument.

\begin{proposition}\label{hyperelliptic}
MRC holds for hyperelliptic curves.
\end{proposition}

\begin{proof}
We show that for every $i$,
$h^{0}(\wedge^{i}Q\otimes \xi)=0$, if $\xi\in {\rm Pic}^{g-2i-1}(C)$
is general.
Since $C$ is hyperelliptic, we have a degree two morphism
$f:C\rightarrow\PP^1$ and if $L=f^*(\OO_{\PP^1}(1))$,
then $\omega_C=L^{g-1}$. Therefore the morphism
$\widetilde{f}\,:\,C\longrightarrow\PP^{g-1}$ defined by $\omega_C$
is the composition of the Veronese embedding $\PP^1\hookrightarrow
\PP^{g-1}$ with $f$. Note that we have $M=\widetilde{f}^*(\Omega_{\PP^{g-1}}
(1))$.

Since on $\PP^1$ we have the exact sequence:
$$0\longrightarrow
\OO_{\PP^1}(-1)^{\oplus(g-1)}\longrightarrow
 H^0(\OO_{\PP^1}(g-1))\tensor\OO_{\PP^1}
\longrightarrow\OO_{\PP^1}(g-1)\longrightarrow 0,$$
we get $M\iso (L^{-1})^{\oplus(g-1)}$.
Therefore for every $i$, we have 
$$\wedge^{i}Q\iso(L^i)^{\oplus{{g-1}\choose i}}.$$
Now if $\xi\in\Pic^{g-2i-1}(C)$ is general,
then $\xi\tensor L^{i}$ is a general line bundle
of degree $g-1$ and so $h^0(\wedge^{i}Q\tensor \xi)=0$.
\end{proof}

\subsection*{\bf Theta divisors and difference varieties for a fixed curve.}
We noted above that MRC is satisfied for $C$ if and only if
 $\Theta_{\wedge^i Q}$ 
is a divisor. We now identify precisely what the divisor should be,
 assuming that 
this happens. (At the end of the day this will hold for all canonical curves.)
Recall that by general theory, whenever a divisor,
$\Theta_{\wedge^i Q}$ belongs 
to the linear series
$|{{g-1} \choose {i}}\Theta|$, where we slightly abusively  
denote by $\Theta$ a certain theta divisor on ${\rm Pic}^{g-2i-1}(C)$ (more
precisely $\Theta_N$, where $N$ is a ${{g-1} \choose {i}}$-th root of 
${\rm det}(\wedge^i Q)$).

From now on we always assume that we are in this situation.
The Picard variety
 ${\rm Pic}^{g-2i-1}(C)$ contains a \emph{difference subvariety}
$C_{g-i-1}-C_{i}$ defined as the image of the difference map
$$\phi:C_{g-i-1}\times C_{i}\longrightarrow {\rm Pic}^{g-2i-1}(C)$$
$$(x_1+\ldots +x_{g-i-1},y_1+\ldots +y_i)\rightarrow \mathcal{O}_C
(x_1+\ldots +x_{g-i-1} -y_1 -\ldots - y_i).$$
The geometry of the difference varieties has interesting links with the 
geometry of the curve \cite{Rob} and \cite{ACGH} (see below). 
The key observation is that our theta divisor 
is nothing else but the difference variety above.

\begin{proposition}\label{1}
For every smooth curve $C$ of genus $g$, we have
$$C_{g-i-1} - C_{i}\subseteq\Theta_{\wedge^i Q}.$$
Moreover, if $C$ is nonhyperelliptic and $\Theta_{\wedge^i Q}$
is a divisor, then the above inclusion is an equality.
\end{proposition}

We start with a few properties of the difference varieties, which for 
instance easily imply that $C_{g-i-1} - C_{i}$ is a divisor. More 
generally, we study the difference variety $C_a - C_b$, $a\geq b$, defined 
analogously. Note that this study is suggested in a series of exercises 
in [ACGH] Ch.V.D and Ch.VI.A in the case $a=b$, but the formula in V.D-3 
there giving the cohomology class of $C_a - C_a$ is unfortunately incorrect,
 as we first learned from R. Lazarsfeld. 
The results we need are collected in the following:

\begin{proposition}\label{2}
(a) Assume that $1\leq b\leq a\leq \frac{g-1}{2}$. Then the difference map:
$$\phi: C_a \times C_b \longrightarrow C_a - C_b\subseteq{\rm
Pic}^{a-b}(C)$$
is birational onto its image if $C$ is
nonhyperelliptic. When $C$ is hyperelliptic, $\phi$ has
 degree ${{a} \choose {b}}\,2^b$ onto its image.
\newline
\noindent
(b) If $C$ is nonhyperelliptic, the cohomology 
class $c_{a,b}$ of $C_a - C_b$ in ${\rm Pic}^{a-b}(C)$ is given by
$$c_{a,b} = {{a+b} \choose {a}}\,\theta^{g-a-b},$$
where $\theta$ is the class of a theta divisor.
\end{proposition}

Assuming this, the particular case $a = g-i-1$ and $b = i$ quickly
implies
the main result.

\begin{proof}(of Proposition \ref{1})
From Proposition \ref{2} (b) we see that if $C$ is nonhyperelliptic, then
the class of $C_{g-i-1} - C_{i}$ is given by:
$$c_{g-i-1,i} = {{g-1} \choose {i}} \theta.$$
On the other hand, as $\Theta_{\wedge^i Q}$ is associated to the
vector bundle $\wedge^i Q$, if it is a divisor, then its cohomology class is 
${{g-1} \choose {i}}\theta$ (recall that $\Theta_{\wedge^i Q}$ has the same 
class as ${{g-1}
\choose {i}}\Theta$). As in this case both $\Theta_{\wedge^i Q}$ and
 $C_{g-i-1} - C_{i}$ are
divisors, in order to finish the proof of the proposition it is enough
to prove the first statement.
 
To this end, we follow almost verbatim the argument in
 Theorem \ref{curves_of_large_degree}
(b). Namely, the filtration (2) in \S2 implies that for every $i\geq 1$
 there is an inclusion:
$$\OO_C(x_1+\ldots + x_i)\hookrightarrow \wedge^i Q,$$
where $x_1,\ldots , x_i$ are general points on $C$. This means that 
$$h^0(\wedge^i Q\otimes \OO_C(E_{g-i-1}-D_i))\neq 0$$
for all general effective divisors $E_{g-i-1}$
 of degree $g-i-1$ and $D_i$ of degree $i$,
which gives the desired inclusion.
\end{proof}

We are left with proving Proposition \ref{2}. This follows by more or less
standard arguments in the study of Abel maps and Poincar\' e formulas for
cohomology classes of images of symmetric products.

\begin{proof}(of Proposition \ref{2})
(a) This is certainly well known (cf. \cite{ACGH} Ch.V.D), and we do not reproduce 
the proof here.
\newline
\noindent
(b) Assume now that $C$ is nonhyperelliptic, so that 
$$\phi: C_a \times C_b \longrightarrow C_a - C_b$$
is birational onto its image. For simplicity we will map 
everything to the Jacobian of $C$, so fix a point $p_0 \in C$ and consider 
the commutative diagram:
$$\xymatrix{
& C^{a+b} \ar[dl]_{\psi} \ar[dr]^{\alpha} \\
C_a\times C_b \ar[r]^{\phi} & C_a - C_b \ar[r]^{-(a-b)p_0} & J(C)}$$
where $C^{a+b}$ is the $(a+b)$-th cartesian product of the curve and the
maps 
are either previously defined or obvious.
We will in general denote by $[X]$ the fundamental class of the compact
variety $X$. Since $\psi$ clearly has degree $a!\cdot b!$, and since
$\phi$ is birational by (a), we have:
$$\alpha_{*}[C^{a+b}] = a!\cdot b!\cdot c_{a,b}.$$
This means that it is in fact enough to prove that $\alpha_{*}[C^{a+b}]
= (a+b)!\cdot \theta^{g-a-b}$, and note that 
$(a+b)!\cdot\theta^{g-a-b}$ is the same as the class
$u_{*}[C^{a+b}]$, where $u$ is the usual Abel map:
$$u: C^{a+b}\longrightarrow J(C)$$
$$(x_1,\ldots , x_{a+b})\rightarrow \OO_C(x_1 + \ldots + x_{a+b} - (a+b)p_0).$$
The last statement is known as Poincar\'e's formula (see e.g. [ACGH] I
\S5).
We are now done by the following lemma,  which essentially says that 
adding or subtracting points is the same when computing cohomology
classes. 
\end{proof}

\begin{lemma}
Let $u,\alpha : C^{a+b} \longrightarrow J(C)$ defined by:
$$u(x_1 ,\ldots ,x_{a+b}) = \OO_C(x_1 + \ldots + x_{a+b} - (a+b)p_0)$$
and 
$$\alpha(x_1 ,\ldots , x_{a+b}) = \OO_C(x_1 + \ldots + x_a - x_{a+1} - \ldots -
x_{a+b} - (a-b)p_0).$$
Then $u_{*}[C^{a+b}] = \alpha_{*}[C^{a+b}]\in H^{2(g-a-b)}(J(C),
\mathbb{Z})$.
\end{lemma}
\begin{proof}
For simplicity, in this proof only, we will use additive divisor notation,
although we actually mean the associated line bundles.
Consider the auxiliary maps:
$$u_0: C^{a+b}\longrightarrow J(C)^{a+b}$$
$$(x_1, \ldots , x_{a+b})\rightarrow (x_1 - p_0 ,\ldots , x_{a+b} -
p_0),$$
$$\alpha_0: C^{a+b}\longrightarrow J(C)^{a+b}$$
$$(x_1,\ldots , x_{a+b})\rightarrow (x_1  - p_0, \ldots , x_a - p_0, 
p_0 - x_{a+1}, \ldots , p_0 - x_{a+b})$$
and the addition map:
$$a: J(C)^{a+b}\longrightarrow J(C)$$
$$(\xi_1, \ldots , \xi_{a+b})\rightarrow \xi_1 + \ldots + \xi_{a+b}.$$
Then one has:
$$u = au_0 ~{\rm and}~\alpha = a\alpha_0 = a\mu u_0,$$
where $\mu$ is the isomorphism
$$\mu: J(C)^{a+b}\longrightarrow J(C)^{a+b}$$
$$(\xi_1, \ldots ,\xi_{a+b})\rightarrow (\xi_1,\ldots,\xi_a , -\xi_{a+1},
\ldots , -\xi_{a+b}).$$
Now the statement follows from the more general fact that $\mu$ induces 
an isomorphism on cohomology. This is in turn a simple consequence
of the fact that the involution $x \rightarrow -x$ induces the identity on
$H^{1}(J(C), \mathbb{Z})$.
\end{proof}

\begin{remark}\label{hyperelliptic_case}
The equality in Proposition~\ref{1} holds set-theoretically also for
hyperelliptic curves. Indeed, we have  
the inclusion $C_{g-i-1}-C_i\subseteq\Theta_{\wedge^iQ}$,
and we have seen in the proof of
Proposition~\ref{hyperelliptic} that $\Theta_{\wedge^iQ}$ is irreducible
(with the reduced structure, it is just a translate of the usual theta
 divisor).
\end{remark}

\subsection*{\bf General canonical curves.}
Since we have seen that MRC holds for hyperelliptic curves, a standard
argument shows that it holds for general canonical curves. In fact, our
previous result about the expected
form of the theta divisors $\Theta_{\wedge^iQ}$ allows us
to say something more precise about the set of curves in ${\mathcal M}_g$
which might not satisfy MRC.

\begin{proposition}\label{generic}
For every $i$, the set of curves
$\{[C]\in {\mathcal M}_g\mid\Theta_{\wedge^iQ_C}=\Pic^{g-2i-1}(C)\}$
is either empty or has pure codimension one. In particular,
the same assertion is true for the set of curves
 in ${\mathcal M}_g$
which do not satisfy MRC.
\end{proposition}

\begin{proof}
As the arguments involved are standard we will just sketch the proof.

Start by considering, for a given $d\geq 2g+1$,   
the Hilbert scheme $\mathcal H$ of curves in $\PP^{d-g}$ with  
Hilbert polynomial $P(T)=d\,T+1-g$ and ${\mathcal U}\subseteq\mathcal{H}$
the open subset corresponding to smooth connected nondegenerate curves.
 
Let $f\,:\,{\mathcal Z}\longrightarrow {\mathcal U}$
 be the universal family over $\mathcal U$,
which is smooth of relative dimension $1$,
 and $\omega_{{\mathcal Z}/{\mathcal U}}
\in\Pic({\mathcal Z})$ the relative cotangent bundle.
By base change there is an exact sequence
$$0\longrightarrow {\mathcal Q}^{\vee}
\longrightarrow f^*f_*\omega_{{\mathcal Z}/
{\mathcal U}}
\longrightarrow 
\omega_{{\mathcal Z}/{\mathcal U}}\longrightarrow 0,$$
where ${\mathcal Q}$ is a vector bundle on $\cZ$ such that 
 if $u\in {\mathcal U}$ corresponds to a curve
 $C=\cZ_u$ (in a suitable embedding), then ${\mathcal Q}\vert_{{\mathcal Z}_u}
\iso Q_C$.

The usual deformation theory arguments show that ${\mathcal H}$ is smooth
and has dimension $(d-g+1)^2+4(g-1)$. Moreover, the universal family
${\mathcal Z}$ defines a surjective
morphism $\pi\,:\,{\mathcal U}\longrightarrow
\cM_g$ whose fibers are irreducible and have dimension $(d-g+1)^2+g-1$.
It is immediate to see from this that $\cU$ is connected.

Fix now $l$ such that $d=(g-2i-1)+l(2g-2)\geq 2g+1$. Consider
$\cU$ and $\cZ$ as above and let $\cF:=\wedge^i\cQ\otimes\omega_{\cZ/\cU}^{-l}
\otimes p^*\cO_{\PP^{d-g}}(1)$, where $p$ is the composition of the
inclusion $\cZ\hookrightarrow\cU\times\PP^{d-g}$ and the projection onto the
second factor.

We consider also the closed subset of $\cU$:
$${\mathcal D}_1=\{u\in\cU\mid h^0(\cF\vert_{\cZ_u})\geq 1\}.$$
It is clear by definition that $\pi^{-1}([C])\subseteq{\mathcal D}_1$
if and only if $\Theta_{\wedge^iQ_C}=\Pic^{g-2i-1}(C)$.
 In particular,
 Proposition~\ref{hyperelliptic} implies that ${\mathcal D}_1\neq \cU$.

${\mathcal D}_1$ is the degeneracy locus of a morphism between two vector
 bundles of the same rank. Indeed, if $H\subset\PP^{d-g}$ is a hyperplane,
$\widetilde{H}=p^{-1}H$,
and $r\gg 0$, then ${\mathcal D}_1$ is the degeneracy locus of 
$$f_*(\cF\otimes\cO_{\mathcal Z}(r\widetilde{H}))
\longrightarrow f_*(\cF\otimes \cO_{r\widetilde{H}}(r\widetilde{H})).$$
Note that these are both vector bundles of rank $rd{{g-1}\choose i}$
(we use base change and the fact that by Corollary~3.5 in 
\cite{paranjape}, for every
smooth curve $C$, the bundle $Q_C$ is semistable). We therefore conclude that
${\mathcal D}_1$ is a divisor on $\cU$.

On the other hand, it is easy to see that the set
$${\mathcal D}_2=\{u\in\cU\mid \cO_{\cZ_u}(1)\otimes\omega_{\cZ_u}^{-l}
\in (\cZ_u)_{g-i-1}-(\cZ_u)_i\}$$ is closed.
Moreover, Proposition~\ref{1}
(see also Remark~\ref{hyperelliptic_case}) shows that
${\mathcal D}_2\subseteq {\mathcal D}_1$, 
and if $\pi^{-1}([C])\not\subseteq{\mathcal D}_1$,
then $\pi^{-1}([C])\cap {\mathcal D}_1=\pi^{-1}([C])\cap {\mathcal D}_2$.

Let $\mathcal S$ be the set of irreducible components of ${\mathcal D}_1$
which are not included in ${\mathcal D}_2$. Using the fact that $\pi$
has irreducible fibers, all of the same dimension, it is easy to see that
if $Y\in{\mathcal S}$, then $\pi(Y)$ is closed in $\cM_g$,
that it is in fact a divisor, and $Y=\pi^{-1}(\pi(Y))$. Moreover, the locus of
curves in $\cM_g$ for which $\Theta_{\wedge^iQ}$ is not a divisor is
$\cup_{Y\in{\mathcal S}}
\pi(Y)$, which proves the proposition.
\end{proof}

\subsection*{\bf The class of the degeneracy locus on $\cM_{g,g+1}$}

We first fix the notation. We will denote by ${\mathcal M}_g^{0}$ the open
subset of ${\mathcal M}_g$ which corresponds to curves with no
nontrivial automorphisms.
From now on we assume that $g\geq 4$, since for $g=3$ MRC is equivalent
to IGC, which is the content of Proposition \ref{igc_canonical}.
Thus $\cM_g^{0}$ is 
nonempty and its complement has codimension 
$g-2\geq 2$ (see \cite{HaMo} pg.37), so working with this subset 
will not affect the answers we get for divisor class computations
on $\cM_g$ or $\cM_{g,n}$.

In this case we have a universal family over $\cM_g^{0}$ denoted by
${\mathcal C}_g^{0}$
 and for every $n\geq 1$, the open subset of
$\cM_{g,n}$ lying over $\cM_g^{0}$
(which we denote by $\cM_{g,n}^{0}$) is equal to the fiber product
$(\times_{\cM_g^{0}}{\mathcal C}_g^{0})^n$
 minus all the diagonals. We assume that $n\geq g+1$.

Consider the following cartezian diagram:
\[
\begin{CD}
{\mathcal X}@>{q}>> \cM_{g,n}^{0}\\
@VV{f}V@VV{h}V\\
{\mathcal C}_g^{0}@>{p}>>\cM_g^{0}\\
\end{CD}
\]
in which all the morphisms 
and varieties are smooth and $p$ (hence also $q$) is proper.

Let $\omega\in\Pic({\mathcal C}_g^{0})$ be the relative canonical line
 bundle for $p$, $E=p_*(\omega)$ the Hodge vector bundle and
 $\cQ$ the rank $g-1$
vector bundle on ${\mathcal C}^{0}_g$ such that $\cQ^{\vee}$ is the
kernel of the evaluation map $p^*E\longrightarrow\omega$. For every
$[C]\in\cM_g^{0}$, we have $\cQ\vert_{p^{-1}([C])}\simeq Q_C$.

The projection on the $j$-th factor $p_j\,:\,\cM_{g,n}^{0}
\longrightarrow{\mathcal C}_g^{0}$ induces a section 
$q_j\,:\,\cM_{g,n}^{0}\longrightarrow {\mathcal X}$ of $q$.
If $E_j={\rm Im}(q_j)$, then $E_j$ is a relative divisor over $\cM_{g,n}
^{0}$. Consider the following vector bundle on ${\mathcal X}$:
$$\cE=\wedge^if^*\cQ\left(\sum_{j=1}^{g-i}E_j-\sum_{j=g-i+1}^{g+1}E_j\right)$$
and let
 $\cZ=\{u\in\cM_{g,n}^{0}\vert ~h^0(\cE\vert_{{\mathcal X}_u})\geq 1\}$.

The algebraic set
$\cZ$ comes equipped with a natural stucture of degeneracy locus. 
Suppose for example that $Y$ is a sum of $m$ divisors $E_j$
(possibly with repetitions), where $m\gg 0$. In this case, $\cZ$
is the degeneracy locus of the morphism
 $$\cF:=q_*(\cE\otimes\cO_{\mathcal X}(Y))
\longrightarrow \cF':=q_*(\cE\otimes\cO_{\mathcal X}(Y)\vert_Y).$$
 This scheme 
structure does not depend on the divisor $Y$ we have chosen.
In fact, it is the universal subscheme over which the $0$-th Fitting ideal
of the first higher direct image of $\cE$ becomes trivial (see e.g.
\cite{ACGH} Ch.IV \S3
 for the proof of an analogous property). Note that $\cZ$ is a divisor
and not the whole space, since by Proposition~\ref{generic} we know that
for a general curve $C$, there is $\xi\in\Pic^{g-2i-1}(C)$
such that $H^0(\wedge^iQ_C\otimes\xi)=0$ (note also 
 that the difference map
$C^{g-i}\times C^{i+1}\longrightarrow \Pic^{g-2i-1}(C)$ is surjective,
cf. \cite{ACGH} Ch.V.D).

We will use the notation $\lambda=c_1(h^*(E))$, $\psi_j=c_1(p_j^{*}(\omega))$,
and $\Psi_x=\sum_{j=1}^{g-i}\psi_j$ and $\Psi_y=\sum_{j=g-i+1}^{g+1}\psi_j$.
It is well known that $\lambda$ together with $\psi_j$, $1\leq j\leq n$,
form a basis for $\Pic(\cM_{g,n}^{0})_{\QQ}$.

\begin{proposition}\label{degeneracy_class}
With the above notation, for every $n\geq g+1$, the class of $\cZ$
in $\Pic(\cM_{g,n}^{0})_{\QQ}$ is given by
\begin{equation}
[\cZ]=-\left({{g-1}\choose i}-10{{g-3}\choose{i-1}}\right)
\lambda+{{g-2}\choose i}\Psi_x
+{{g-2}\choose{i-1}}\Psi_y.
\end{equation}
\end{proposition}

\begin{proof}
Note that the pull-back of divisors induced by the projection to
 the first $(g+1)$ components induces injective homomorphisms
$\Pic(\cM_{g,g+1}^{0})_{\QQ}\hookrightarrow
\Pic(\cM_{g,n}^{0})_{\QQ}$. From the universality of the scheme
structure on $\cZ$ it follows that
 the computation of $c_1(\cZ)$ is
independent of $n$. Therefore we may assume
that $n$ is large enough, so that in defining the
scheme structure of $\cZ$ as above, we may take $Y=\sum_{j=g+2}^nE_j$.
We introduce also the notation $\Psi_z=\sum_{j=g+2}^n\psi_j$.

As a degeneracy locus, the class of $\cZ$ is given by $c_1(\cF')-c_1(\cF)$.
It is clear that we have $E_j\cap E_l=\emptyset$ if $j\neq l$
and via $E_j\simeq\cM_{g,n}^{0}$, we have $\cO_{E_j}(-E_j)\simeq
p_j^*(\omega)$. Since $\cQ^{\vee}$ is the kernel of the evaluation map for
$\omega$, we get 
\begin{equation}\label{first_formula}
c_1(f^*Q)=f^*(c_1(\omega))-q^*(\lambda).
\end{equation}

Before starting the computation of $c_1(\cF)$ and $c_1(\cF')$, we record the
following well-known formulas for Chern classes.

 \begin{lemma}\label{chern_classes}
Let $R$ be a vector bundle of rank $n$ on a variety $X$
and $L\in {\rm Pic}(X)$. 
\item{\rm (i)} $c_1(R\otimes L)=c_1(R)+nc_1(L)$.
\item{\rm (ii)} $c_2(R\otimes L)=
c_2(R)+(n-1)c_1(R)c_1(L)+{n\choose 2}c_1(L)^2$.
\item{\rm (iii)} $c_1(\wedge^iR)=
{{n-1}\choose{i-1}}c_1(R)$, if $n\geq 2$ and $1\leq i
\leq n$.
\item{\rm (iv)} $c_2(\wedge^iR)=\frac12
{{n-1}\choose{i-1}}\left({{n-1}\choose {i-1}}-1\right)
c_1(R)^2+{{n-2}\choose{i-1}}c_2(R)$, if $n\geq 3$ and $1\leq i\leq n$.
\end{lemma}

{}From the previous discussion and Lemma~\ref{chern_classes} (i)
 and (iii), we get
$$c_1(q_*(\cE\otimes\cO_{E_j}(Y)))=-{{g-2}\choose{i-1}}\lambda-
{{g-2}\choose i}\psi_j,$$
{}from which we deduce
$$c_1(\cF')=-(n-g-1){{g-2}\choose{i-1}}\lambda-{{g-2}\choose i}\Psi_z.$$

In order to compute $c_1(\cF)$, we apply the Grothendieck-Riemann-Roch
formula for $q$ and $\cE(Y)$
(see \cite{fulton} 15.2).
 Note that the varieties are smooth and $q$ is smooth and proper.
Moreover, since we assume $n\gg 0$, we have $R^jq_*(\cE(Y))=0$, for $j\geq 1$.
Therefore we get
$$ch(q_*(\cE(Y)))=q_*(ch(\cE(Y))\cdot td(f^*\omega^{-1})).$$
From this we deduce
$$(*)\,
c_1(\cF)=q_*\left(\frac12 c_1(\cE(Y))^2-c_2(\cE(Y))-
\frac12 f^*c_1(\omega)\cdot c_1(\cE(Y))
+\frac1{12}{{g-1}\choose i}f^*c_1(\omega)^2\right). 
$$

We compute now each of the classes involved in the above equation.
In order to do this we need to know how to make the push-forward of the
elementary classes on $\mathcal X$. We list these rules in the following: 

\begin{lemma}\label{rules}
With the above notation, we have
\item {\rm (i)} $q_*(f^*c_1(\omega)^2)=12\lambda$.
\item {\rm (ii)} $q_*(q^*\lambda\cdot f^*c_1(\omega))=(2g-2)\lambda$.
\item {\rm (iii)} $q_*q^*(\lambda^2)=0$.
\item {\rm (iv)} $q_*(c_1(E_j)\cdot q^*\lambda)=\lambda$.
\item {\rm (v)} $q_*(c_1(E_j)\cdot f^*c_1(\lambda))=\psi_j$.
\item {\rm (vi)} $q_*q^*c_2(h^*E)=0$.
\item {\rm (vii)} $q_*(c_1(E_j)^2)=-\psi_j$.
\end{lemma}

\begin{proof}[Proof of Lemma~\ref{rules}]
The proof of (i) is analogous to that of the relation $p_*(c_1(\omega)^2)=
12\lambda$ (see \cite{HaMo} \S 3E). The other formulas are straightforward.
\end{proof}

Using Lemma~\ref{chern_classes} (i) and (iii) and the formula 
(\ref{first_formula}) for
$c_1(f^*Q)$, we deduce that 
$$c_1(\cE(Y))={{g-2}\choose{i-1}}(f^*c_1(\omega)-q^*\lambda)
+{{g-1}\choose i}\left(\sum_{j=1}^{g-i}c_1(E_j)-\sum_{j=g-i+1}^{g+1}c_1(E_j)
+\sum_{j=g+2}^nc_1(E_j)\right).$$
Applying Lemma~\ref{rules}, we get
$$q_*(c_1(\cE(Y))^2/2)=\left((8-2g){{g-2}\choose{i-1}}^2
-(n-2i-2){{g-2}\choose{i-1}}
{{g-1}\choose i}\right)\lambda+$$
$${{g-2}\choose{i-1}}{{g-1}\choose i}(\Psi_x-\Psi_y+\Psi_z)
-\frac12{{g-1}\choose i}^2(\Psi_x+\Psi_y+\Psi_z).$$

{}From the above formula for $c_1(\cE(Y))$ and Lemma~\ref{rules}, we get
$$q_*\left(-\frac12 f^*c_1(\omega)\cdot c_1(\cE(Y))\right)
=(g-7){{g-2}\choose{i-1}}\lambda
-\frac12{{g-1}\choose i}(\Psi_x-\Psi_y+\Psi_z).$$

Lemma~\ref{rules} (i) gives
$$q_*\left(\frac1{12}{{g-1}\choose i}f^*c_1(\omega)^2\right)
={{g-1}\choose i}\lambda.$$

{}From the defining exact sequence of $\cQ^{\vee}$ we compute
$$c_2(f^*\cQ)=q^*c_2(h^*E)+f^*c_1(\omega)\cdot(f^*c_1(\omega)-q^*\lambda).$$
Using now Lemma~\ref{chern_classes} (ii) and (iv) and Lemma~\ref{rules},
 we deduce
$$q_*c_2(\cE(Y))= \left((8-2g){{g-2}\choose{i-1}}\left(
{{g-2}\choose{i-1}}-1\right)+(14-2g){{g-3}\choose{i-1}}\right)\lambda-$$
$$(n-2i-2){{g-2}\choose{i-1}}\left({{g-1}\choose i}-1\right)\lambda+
{{g-2}\choose{i-1}}\left({{g-1}\choose i}-1\right)(\Psi_x-\Psi_y+\Psi_z)-$$
$$\frac12{{g-1}\choose i}\left({{g-1}\choose i}-1\right)(\Psi_x
+\Psi_y+\Psi_z).$$

Using these formulas and equation $(*)$, we finally obtain
$$c_1(\cF)=\left((2g-14){{g-3}\choose{i-1}}-(n+g-2i-3){{g-2}\choose{i-1}}
+{{g-1}\choose i}\right)\lambda-$$
$${{g-2}\choose i}(\Psi_x+\Psi_z)
-{{g-2}\choose{i-1}}\Psi_y.$$

Since the class of $\cZ$ is equal with $c_1(\cF')-c_1(\cF)$, we deduce
the statement of the proposition.
\end{proof}

\subsection*{\bf Proof of the main result.}
We introduce next a divisor
 $D$ on $\mathcal{M}_{g,g+1}$ which is a global analogue of the
 preimage of $C_{g-i-1}-C_{i}$ under the difference map
 $C^{g-i}\times C^{i+1}\rightarrow \mbox{Pic}^{g-2i-1}(C)$.
 This is motivated by Proposition \ref{1}, and our goal is roughly 
 speaking to prove a global version of that result. 

\begin{definition} For $g\geq 3$ and
 $1\leq i\leq \frac{g-1}{2}$ we define the divisor $D$
on $\mathcal{M}_{g,g+1}$ to be the locus of smooth pointed curves
 $(C, x_1,\ldots ,x_{g-i}, y_1,\ldots y_{i+1})$
 having a linear series $\mathfrak g^1_{g}$ containing 
$x_1+\cdots +x_{g-i}$ in a fiber and $y_1+\cdots +y_{i+1}$ in another fiber.
Note that this means that we can in fact write the line bundle
$\OO_C(x_1+\ldots +x_{g-i}-y_1-\ldots-y_{i+1})$ 
as an element in $C_{g-i-1}-C_i$.
\end{definition}

We consider in what follows the divisor $\overline{\mathcal{Z}}$,
the closure in ${\mathcal M}_{g,g+1}$ of the divisor 
${\mathcal Z}$ studied above (we take now $n=g+1$,  but
as we mentioned, this does not affect the formula for its class).
In Section 4 we prove that $D$ is reduced and that
 $D\equiv_{\QQ-{lin}} \overline{\mathcal{Z}}$
(cf. Theorem \ref{class_D}).
This being granted we are in a position to complete the proof of Theorem 3.1:

\begin{proof}[Proof of Theorem~\ref{canonical_curves}] 
Note that for $g=3$ our assertion is just the statement of Proposition 
\ref{igc_canonical}. Thus we can assume $g\geq 4$.
As mentioned above $D$ is reduced, and from Proposition \ref{1}
we see that
 $\mbox{supp}(D)\subseteq \mbox{supp}(\overline{\mathcal{Z}})$.
 We get that $\overline{\mathcal{Z}}-D$ is effective,
 and in fact $\overline{\mathcal{Z}}-D=h^*(E)$, where $E$ is an effective
 divisor on $\mathcal{M}_g$ and $h:\mathcal{M}_{g,g+1}\rightarrow
 \mathcal{M}_g$ is the projection. 
Moreover, the map $h^*:\mbox{Pic}(\mathcal{M}_g)_{\QQ}\rightarrow
 \mbox{Pic}(\mathcal{M}_{g,g+1})_{\QQ}$ is injective (cf. \cite{AC2}),
 hence $E\equiv_{\QQ-lin}0$. Since 
the Satake compactification of $\mathcal{M}_g$ has boundary
 of codimension $2$ (see e.g \cite{HaMo}, pg.45) this implies $E=0$,
 that is, $\overline{\mathcal{Z}}=D$.
 Therefore $\Theta_{\wedge ^i Q_C}$ is a divisor in
 $\mbox{Pic}^{g-2i-1}(C)$ and the identification
 $\Theta_{\wedge ^i Q_C}=C_{g-i-1}-C_{i}$ holds for
 \emph{every} nonhyperelliptic curve $C$.
\end{proof}

\section{\bf A divisor class computation on $\Mg1$}

In this section we compute the class of the divisor
 $D$ on $\mathcal{M}_{g,g+1}$ defined in the previous section.
We start by recalling a few facts about line bundles on
 $\overline{\mathcal{M}}_{g,n}$. Let us fix $g\geq 3$, $n\geq 0$
 and a set $N$ of $n$ elements.
 Following \cite{AC2}, we identify $\overline{\mathcal{M}}_{g,n}$
 with the moduli space
 $\overline{\mathcal{M}}_{g,N}$
 of stable curves of genus $g$ with marked points indexed by $N$.
 We denote by $\pi_{q}:\overline{\mathcal{M}}_{g,N\cup\{q\}}\rightarrow
 \overline{\mathcal{M}}_{g,N}$ the map forgetting the marked
 point indexed by $q$.
 For each $z\in N$ we define the tautological class $\psi_z=c_1(\mathbb L_z)\in \mbox{Pic}(\overline{\mathcal{M}}_{g,N})_{\mathbb Q}$,
where $\mathbb L_z$ is the line bundle over $\overline{\mathcal{M}}_{g,N}$ whose fiber over the moduli point $[C,\{x_i\}_{i\in N}]$ is the cotangent space $T_{x_z}^*(C)$. Note that although we are
 using an apparently different definition, these $\psi$ classes are the same 
 as those which appear in the previous section.

 For $0\leq i\leq g$ and $S\subseteq N$, the boundary divisor
 $\Delta_{i:S}$  corresponds to the closure in $\overline{\mathcal{M}}_{g,N}$
of the locus of nodal curves
 $C_1\cup C_2$, with $C_1$ smooth of genus $i$, $C_2$ smooth
 of genus $g-i$, and such that the marked points sitting on
 $C_1$ are precisely those labelled by $S$.
 Of course $\Delta_{i:S}=\Delta_{g-i:S^{c}}$ and we set
 $\Delta_{0:S}:=0$ when $|S|\leq 1$. We also consider the divisor $\Delta_{irr}$ consisting of irreducible pointed curves
with one node. We denote by $\delta_{i:S}\in \mbox{Pic}(\overline{\mathcal{M}}_{g,n})_{\QQ}$ the class of $\Delta_{i:S}$ and by $\delta_{irr}$ that of $\Delta_{irr}$. 
 It is well known that the Hodge class $\lambda, \delta_{irr}$, the $\psi_z$'s and the
 $\delta_{i:S}$'s freely generate
 $\mbox{Pic}(\overline{\mathcal{M}}_{g,n})_{\QQ}$ (cf. \cite{AC2}). 

 For a smooth curve $C$ and for a pencil $\mathfrak g^1_d$ on $C$, we say that an effective divisor $E$ on $C$ is in a fiber of the pencil if there exists 
$E'\in \mathfrak g^1_d$ such that $E'-E$ is an effective divisor.

Recall that for $g\geq 3$ and
 $0\leq i\leq \frac{g-1}{2}$ we have defined the divisor $D$
on $\mathcal{M}_{g,g+1}$ to be the locus of curves
 $(C, x_1,\ldots ,x_{g-i}, y_1,\ldots ,y_{i+1})$
 having a linear series $\mathfrak g^1_{g}$ containing 
$x_1+\cdots +x_{g-i}$ in a fiber and $y_1+\cdots +y_{i+1}$ in another fiber.
We denote by $\overline{D}$ the closure of $D$ in
 $\overline{\mathcal{M}}_{g,g+1}$.

The divisor $D$ comes equipped with a scheme structure induced by the forgetful map $\mathcal{G}\rightarrow \mathcal{M}_{g,g+1}$. Here $\mathcal{G}$ is the variety parametrizing objects $[C,\vec{x},\vec{y},l]$, where $\vec{x}=(x_1,\ldots,x_{g-i})$ and $\vec{y}=(y_1,\ldots, y_{i+1})$ are such that $[C,\vec{x},\vec{y}]\in \mathcal{M}_{g,g+1}$ and $l$ is a linear series $\mathfrak g^1_g$ on $C$ such that $l(-\sum_{j=1}^{g-i} x_j)\neq \emptyset$ and
$l(-\sum_{j=1}^{i+1} y_j)\neq \emptyset$.
It is well-known that $\mathcal{G}$ is smooth of pure dimension $4g-3$ (see e.g. \cite{AC1}, pg. 346). Note also that there is a natural action of $S_{g-i}\times S_{i+1}$  on $\mathcal{M}_{g,g+1}$ (and hence on $\mathcal{G}$) by permuting the components of $\vec{x}$ and $\vec{y}$ separately.

The main result of the section is the following:

\begin{theorem}\label{class_D}
The divisor $D$ is reduced and its class in $\Pic(\overline{\mathcal{M}}_
{g,g+1})_{\QQ}$ is
$$[D]=-\Bigl({g-1\choose  i}
-10{g-3\choose i-1}\Bigr)\ \lambda+{g-2\choose i}\Psi
_x+{g-2\choose i-1}\Psi_y, $$
where $\Psi_x=\sum_{j=0}^{g-i} \psi_{x_j}$
 and $\Psi_y=\sum_{j=0}^{i+1} \psi_{y_j}$.
\end{theorem}

\noindent
We begin by proving the first part of Theorem~\ref{class_D}:

\begin{proposition}\label{reduced}
The divisor $D$ is reduced.
\end{proposition}

 \begin{proof}
 Since the variety $\mathcal{G}$ introduced above is smooth, it suffices to show that the projection $\pi:\mathcal{G}\rightarrow
 \mathcal{M}_{g,g+1}$ given by $\pi([C,\vec{x},\vec{y},l])=[C,\vec{x},\vec{y}]$,
is generically injective. 

We pick a component $\mathcal{X}$ of $\mathcal{G}\times _{\mathcal{M}_{g,g+1}}\mathcal{G}$ whose general point corresponds to a marked curve $[C,\vec{x},\vec{y}]\in \mathcal{M}_{g,g+1}$ together with two {\emph{different}} base-point-free $\mathfrak g^1_g$'s on $C$, both containing $\vec{x}$ and $\vec{y}$ in different fibers. Clearly $\mbox{dim}(\mathcal{X})\geq 4g-4$ and if we show that $\mbox{dim}(\mathcal{X})\leq 4g-4$, then we are done.
For a general point in $\mathcal{X}$ we denote by $f_1,f_2:C\rightarrow \PP^1$ the induced $g$-sheeted maps. We may assume that $f_1(\vec{x})=f_2(\vec{x})=0$ and $f_1(\vec{y})=f_2(\vec{y})=\infty.$
The product map $f=(f_1,f_2):C
\rightarrow \PP^1\times \PP^1$ is birational onto its image
 and $\Gamma=f(C)$ will have points of multiplicity
 at least $g-i$ and $i+1$ at $a=(0,0)$ and $b=(\infty,\infty)$ respectively.

 If $S=\mbox{Bl}_{\{a,b\}}(\PP^1\times \PP^1)$ we set $\gamma= gl+gm-(g-i)E_a-(i+1)E_
b\in \mbox{Pic(S)}$, where $l$ and $m$ are pullbacks of
 the rulings on $\PP^1\times \PP^1$,
 and $E_a$, $E_b$ are the exceptional divisors. We denote by $V(S,\gamma)$ 
the Severi variety of curves $Y\subset S$ homologous to $\gamma$. 
The discussion above shows that $\mathcal{X}$ lies in the closure of the image of the rational map $V(S,\gamma)- ->(\mathcal{G}\times_{\mathcal{M}_{g,g+1}}\mathcal{G})/S_{g-i}\times S_{i+1}$
obtained by projecting $S$ onto the two factors. Thus $\mbox{dim}(\mathcal{X})\leq \mbox{dim }V(S,\gamma)-\mbox{dim} \mbox{ Aut}(S).$

On the other hand an argument identical to that in \cite{AC1}, Proposition 2.4, 
shows that since $S$ is a regular surface, every irreducible component $M$ of $V(S,\gamma)$ 
having dimension $\geq g+1$ is of the expected dimension provided by deformation theory, 
that is, 
$\mbox{dim}(M)=g-1-\gamma \cdot K_S$. 
Therefore $\mbox{dim}(\mathcal{X})\leq g-1-\gamma \cdot K_S-
\mbox{dim }\mbox{Aut}(S)=4g-4$.
\end{proof}

We will prove the second part of Theorem~\ref{class_D} using degeneration
techniques and enumerative geometry.

\subsection*{Recap on limit linear series (cf. \cite{EH1})}
We recall that for a smooth curve $C$,
 a point $p\in C$ and a linear series $l=(L,V)$
 with $L\in \mbox{Pic}^d(C)$ and $V\in G(r+1, H^0(L))$,
 the \emph {vanishing sequence} of $l$ at $p$
 is obtained by ordering the set $\{\mbox{ord}_p(\sigma)\}_{\sigma \in
 V}$, and it is denoted by 
$$a^l(p):0\leq a_0^l(p)<\ldots <a_r^l(p)\leq d.$$
The \emph{weight} of $l$ at $p$ is defined as
 $w^l(p):=\sum_{i=0}^r (a_i^l(p)-i)$. 

 Given a curve $C$ of compact type, a \emph{limit}
 $\mathfrak g^r_d$ on $C$ is a collection of honest
 linear series $l_Y=(L_Y,V_Y)\in G^r_d(Y)$ for each
 component $Y$ of $C$, satisfying the compatibility
 condition that if $Y$ and $Z$ are components of $C$ 
meeting at $p$ then
$$a_i^{l_Y}(p)+a_{r-i}^{l_Z}(p)\geq d\mbox{ }\mbox{ for }i=0,\ldots ,r.$$
We note that limit linear series appear as limits of
 ordinary linear series in $1$-dimensional families of
 curves and there is a useful sufficient criterion
 for a limit $\mathfrak g^r_d$ to be
 \emph{smoothable} (cf. \cite{EH1}, Theorem 3.4).

\smallskip
\noindent
We will need the following enumerative result (cf. \cite{H}, Theorem 2.1):

\begin{proposition}\label{enumerative}
Let $C$ be a general curve of genus $g$,
 $d\geq \frac{g+2}{2}$ and $p\in C$ a general point.
\begin{itemize}
\item The number of $\mathfrak g^1_d$'s on $C$ containing
 $(2d-g)q$ in a fiber, where $q\in C$ is an unspecified point, is
$$b(d,g)=(2d-g-1)(2d-g)(2d-g+1)\frac{g!}{d!(g-d)!}\mbox{ }.$$
\item If $\beta \geq 1$, $\gamma \geq 1$ are integers
 such that $\beta+\gamma=2d-g$, the number of
 $\mathfrak g^1_d$'s on $C$ containing $\beta p+\gamma q$
 in a fiber for some point $q \in C$ is
$$c(d,g,\gamma)=(\gamma^2(2d-g)-\gamma)\frac{g!}{d!\ (g-d)!}\ .$$
\end{itemize}
\end{proposition}

\noindent
The following simple observation will be used repeatedly:

\begin{proposition}\label{glue}
Fix $y,z\in N$ and denote by
 $\pi_z:\overline{\mathcal{M}}_{g,N}\rightarrow
 \overline{\mathcal{M}}_{g,N-\{z\}}$ the map forgetting
 the marked point labelled by $z$. If $E$ is any divisor class
 on $\overline{\mathcal{M}}_{g,N}$, then the $\lambda$ and the $\psi_x
$ coefficients of $E$ are the same as those of
 $(\pi_z)_*(E\cdot \delta_{0:yz})$ for all $x\in \{y,z\}^c.$
\end{proposition}

\begin{proof}
 We write $E$ uniquely as a combination of $\lambda$, tautological classes 
$\psi_y,\psi_z$ and $\psi_x$ with $x\in \{y,z\}^c$ and boundary divisors. To express $(\pi_z)_*(E\cdot \delta_{0:yz})$ in $\mbox{Pic}(\overline{\mathcal{M}}_{g,N-\{z\}})_{\mathbb Q}$ we use that
$(\pi_z)_*(\lambda\cdot \delta_{0:yz})=
\lambda,\mbox{ }(\pi_z)_*(\psi_x\cdot \delta_{0:yz})
=\psi_x \mbox{ for }x\in \{y,z\}^c$ and that $(\pi_z)_*(\psi_x\cdot \delta_{0:yz})=0$ for $x\in \{y,z\}.$ Moreover we have that $(\pi_z)_*(\delta_{i:S}\cdot \delta_{0:yz})$ is boundary in all cases
except that $(\pi_z)_*(\delta_{0:yz}^2)=-\psi_y$ (cf. \cite{AC2}, Lemma 1.2
 and \cite{logan}, Theorem 2.3). The conclusion follows immediately. 
\end{proof}

 By a succession of push-forwards, using Proposition \ref{glue}
 we will reduce the problem of computing the class of $D$
 to two divisor class computations in $\overline{\mathcal{M}}_{g,3}$. The main idea is to let all the points $x_j$ and then all the points $y_j$ come together and understand how the geometric condition defining $D$ changes under degeneration. Recall that by $\overline{D}$ we denote the closure of $D$ in $\overline{\mathcal{M}}_{g,g+1}$.

 We define the following sequence of divisors:
 starting with $\overline{D}=D_{y_{i+1}}$, for $1\leq j\leq i$
 we define inductively the divisors $D_{y_j}$ on
 $\overline{\mathcal{M}}_{g,g-i+j}$ by
$$D_{y_j}:=(\pi_{y_{j+1}})_*(\Delta_{0:y_jy_{j+1}}\cdot D_{y_{j+1}}).$$
Loosely speaking, $D_{y_j}$ is obtained from $D_{y_{j+1}}$
 by letting the marked points $y_j$ and $y_{j+1}$ come together.
Then we define $D_{x_{g-i}}:=D_{y_1}$ and we let the
 marked points $x_2,\ldots,x_{g-i}$ come together
: for $2\leq j\leq g-i-1$
 we define inductively the divisors $D_{x_j}$ on
 $\overline{\mathcal{M}}_{g,j+1}$ by
$$D_{x_j}:=(\pi_{x_{j+1}})_*(\Delta_{0:x_j x_{j+1}} \cdot D_{x_{j+1}}).$$
Proposition~\ref{glue}
 ensures that the $\psi_{x_1}$ and the $\lambda$ coefficients
of $[\overline{D}]$ are the same as those of $[D_{x_2}]$.

\begin{proposition}\label{reduce}
 The divisor $D_{x_2}$ is reduced and
 it is the closure in $\overline{\mathcal{M}}_{g,3}$
 of the locus of those smooth pointed curves $(C,x_1,x_2,y)$
 for which there exists a $\mathfrak g^1_g$ with
 $(i+1)y$ in a fiber and $x_1+(g-i-1)x_2$ in another fiber.
\end{proposition}

\begin{proof}
 For simplicity we will only prove that
 $D_{y_i}$ is reduced and that it is the closure of the locus
 of those smooth pointed curves $(C,x_1,\ldots,x_{g-i},y_1,\ldots,y_i)$
 for which $x
_1+\cdots +x_{g-i}$ and $y_1+\cdots +y_{i-1}+2y_i$ are
 in different fibers of the same $\mathfrak g^1_g$.
 Then by iteration we will get a similar statement for $D_{x_2}$.

 Let $(X=C\cup_q \PP^1,x_1,\ldots ,x_{g-i}, y_1,\ldots ,y_{i+1})$
 with $y_i, y_{i+1}\in \PP^1$ be a general point in a
 component of $D_{y_{i+1}}\cap \Delta_{0:y_i y_{i+1}}$. A standard dimension count shows that $C$ must be smooth. There exists a limit $\mathfrak g^1_g$ on $X$,
 say $l=(l_C,l_{\PP^1})$, together with sections
 $\sigma_{\PP^1} \in V_{\PP^1}$ and
 $\sigma_C, \tau_C \in V_C$, such that
 $\mbox{div}(\tau_C)\geq x_1+\cdots +x_{g-i}$,
 $\mbox{div}(\sigma_C)\geq y_1+\cdots +y_{i-1}$,
 $\mbox{div}(\sigma_{\PP^1})\geq y_i+y_{i+1}$ and moreover
 $\mbox{ord}_q(\sigma_{\PP^1})+\mbox{ord}_{q}(\sigma_C)\geq g$
 (apply \cite{EH1}, Proposition 2.2). 

Clearly $\mbox{ord}_q(\sigma_{\PP^1})\leq g-2$, hence
 $\mbox{div}(\sigma_C)\geq 2q+y_1+\cdots +y_{i-1}$.
 The contraction map $\pi_{y_{i+1}}$ collapses $\PP^1$
 and identifies $q$ and $y_i$, so the second part of the claim follows.

 To conclude that $D_{y_i}$ is also reduced we use
 that both $D_{y_{i+1}}$ and $\Delta_{0:y_i y_{i+1}}$
 are reduced and that they meet transversally.
 This is because the limit $\mathfrak g^1_g$ we found on $X$
 is smoothable in such a way that all ramification is kept
 away from the nodes (cf. \cite{EH1}, Proposition 3.1), hence the tangent spaces to $D_{y_i}$ and $\Delta_{0:y_iy_{i+1}}$ at the intersection point $(X,\vec{x},\vec{y})$ cannot be equal.
\end{proof}

\smallskip

 In a similar way, by letting first all $x_j$ with
 $1\leq j\leq g-i$ and then all $y_j$ with $2\leq j\leq i+1$
 coalesce, we obtain a reduced divisor $D_{y_2}$ on
 $\overline{\mathcal{M}}_{g,3}$ which is the closure of the locus
 of smooth curves $(C,x,y_1,y_2)$ having a $\mathfrak g^1_g$
 with $(g-i)x$ and $y_1+iy_2$ in different fibers.
 Moreover, the $\lambda$ and the $\psi_{y_1}$ coefficients of
 $[\overline{D}]$ coincide with those of $[D_{y_2}]$.
 Once more applying Proposition~\ref{glue} it follows that the $\lambda$
 and the $\psi_{y_1}$ coefficients of $[D_{y_2}]$ are the same
 as those of 
 $(\pi_x)_*([D_{y_2}]\cdot \delta_{0:xy_2})$. Similarly, 
the $\psi_{x_1}$ coefficient of $[D_{x_2}]$ is the same as that
 of $(\pi_y)_*([D_{x_2}]\cdot \delta_{0:x_2y})$.
 
\begin{proposition}\label{sum_Y}
 We have that 
$$(\pi_x)_*(D_{y_2}\cdot \Delta_{0:xy_2})=\sum_{j=0}^i
Y_j,$$
where for $j<i$ the reduced divisor $Y_j$ is the closure in
 $\overline{\mathcal{M}}_{g,2}$ of the locus of curves
 $(C,y_1,y_2)$ having a $\mathfrak g^1_{g-j}$ with
 $(g-2j-1)y_2+y_1$ in a fiber, while the reduced divisor $Y_i$ consists of curves
 $(C,y_1,y_2)$ with
a $\mathfrak g^1_{g-i}$ having $(g-2i)y_2$ in a fiber
 (and no condition on $y_1$).
\end{proposition}

\begin{proof}
 Once again, let $(X=C\cup_q \PP^1,x,y_1,y_2)$
 be a point in $D_{y_2}\cap \Delta_{0:xy_2}$, with $y_1\in C$
 and $x,y_2\in \PP^1$. Then there exists a limit
 $\mathfrak g^1_g$, say $l=(l_C,l_{\PP^1})$ on $X$ together
 with sections $\sigma_{\PP^1},\tau_{\PP^1}\in
 V_{\PP^1}$ and $\sigma_C\in V_C$ such that 
$\mbox{div}(\sigma_{\PP^1})\geq iy_2, \mbox{div}
(\tau_{\PP^1})\geq (g-i)x, \mbox{div}(\sigma_C)\geq y_1$
 and moreover $\mbox{ord}_q(\sigma_{\PP^1})
+\mbox{ord}_q(\sigma_C)\geq g$.

 The Hurwitz formula on $\PP^1$ and the condition defining 
 a limit linear series give that
 $w^{l_C}(q)\geq w^{l_{\PP^1}}(x)+w^{l_{\PP^1}}(y_2)\geq g-2$.
 On the other hand, since $(X,x,y_1,y_2)$ moves in a family
 of dimension $\geq 3g-2$ it follows that $(C,q)$ also moves in a
 family of dimension $\geq 3g-3$ in $\mathcal{M}_{g,1}$
 (i.e. codimension $\leq 1$). Since according to \cite{EH2}, Theorem 1.2,
 the locus of pointed curves $[C,q]\in \mathcal{M}_{g,1}$ carrying a $\mathfrak g^1_g$ having $w(q)\geq g$ has codimension $\geq 2$, we get $w^{l_C}(q)\leq g-1$. There are two possibilities:

 \textbf{i)} $w^{l_C}(q)=g-2$. Let us denote $j=a_0^{l_C}(q)$,
 hence $a_1^{l_C}(q)=g-1-j$ and $a_k^{l_C}(q)+a_{1-k}^{l_{\PP^1}}(q)=g$
 for $k=0,1$. Therefore
 $j+1=a_0^{l_{\PP^1}}(q)\leq \mbox{ord}_q(\tau_{\PP^1})\leq i.$
 Moreover, since $\mbox{ord}_{q}(\sigma_{\PP^1})
\leq g-i\leq g-j-1$, we obtain that $\mbox{ord}_q(\sigma_C)\geq j+1$,
 hence $\mbox{div}(\sigma_C)\geq y_1+(g-1-j)q$,
 that is, $l_C(-jq)$ is a $\mathfrak g^1_{g-j}$
 on $C$ with $(g-2j-1)q+y_1$ in a fiber, or equivalently $
[C,y_1,q]\in Y_j$, where $0\leq j\leq i-1$.

To see that conversely $\bigcup_{j=0}^{i-1}Y_j\subseteq
 (\pi_{x})_*(D_{y_2}\cdot \Delta_{0:xy_2})$ we pick a general
 pointed curve $(C,y_1,q)$ having a $\mathfrak g^1_{g-j}$ with
 $(g-2j-1)q+y_1$ in a fiber and we construct a Harris-Mumford admissible
 covering $f:X'\rightarrow B$ of degree $g$,
 where $X'$ is a curve semistably equivalent to $X$ defined as above, and
 $B=(\PP^1)_1\cup_{t}(\PP^1)_2$ is the transversal
 union of two lines (see Fig. 1): we take $f_{|C}:C\rightarrow (\PP^1)_1$
 to be the degree $g-j$ covering 
such that $(g-2j-1)q+y_1\subseteq f_{|C}^*(t)$, while
 $f_{|\PP^1}:\PP^1\rightarrow (\PP^1)_2$
 is the degree $g-j-1$ map containing $(g-i)x$ and $iy_2$
 in different fibers and with $(g-2j-1)q$ in the fiber over $t$.
 It is clear that there
 is a unique such $\mathfrak g^1_{g-j-1}$ on $\PP^1$.
 Furthermore, at $y_1$ we insert a rational curve $R$ mapping
 isomorphically onto $(\PP^1)_2$ and at the remaining $j$ points
 in $f_{|C}^{-1}(t)-\{y_1,q\}$ we insert rational curves mapping 
with degree $1$ onto $(\PP^1)_2$ while at the $g-j$ points
 in $f_{|\PP^1}^{-1}(t)-\{q\}$ we insert
 copies of $\PP^1$ mapping isomorphically onto $(\PP^1)_1$. 
We denote the resulting curve by $X'$. If $y_1'=f_{|R}^{-1}(f(y_2))$, then
$(X',x,y_1',y_2)$ is stably equivalent to
 $(X,x,y_1,y_2)$ and $iy_2+y_1'$ and $(g-i)x$ appear in distinct fibers of the $g$-sheeted map $f:X'\rightarrow B$. Thus we get that $[X,x,y_1,y_2]\in D_{y_2}\cap \Delta_{0:xy_2}$.

\begin{figure}[ht]
\begin{center}
\mbox{\epsfig{file=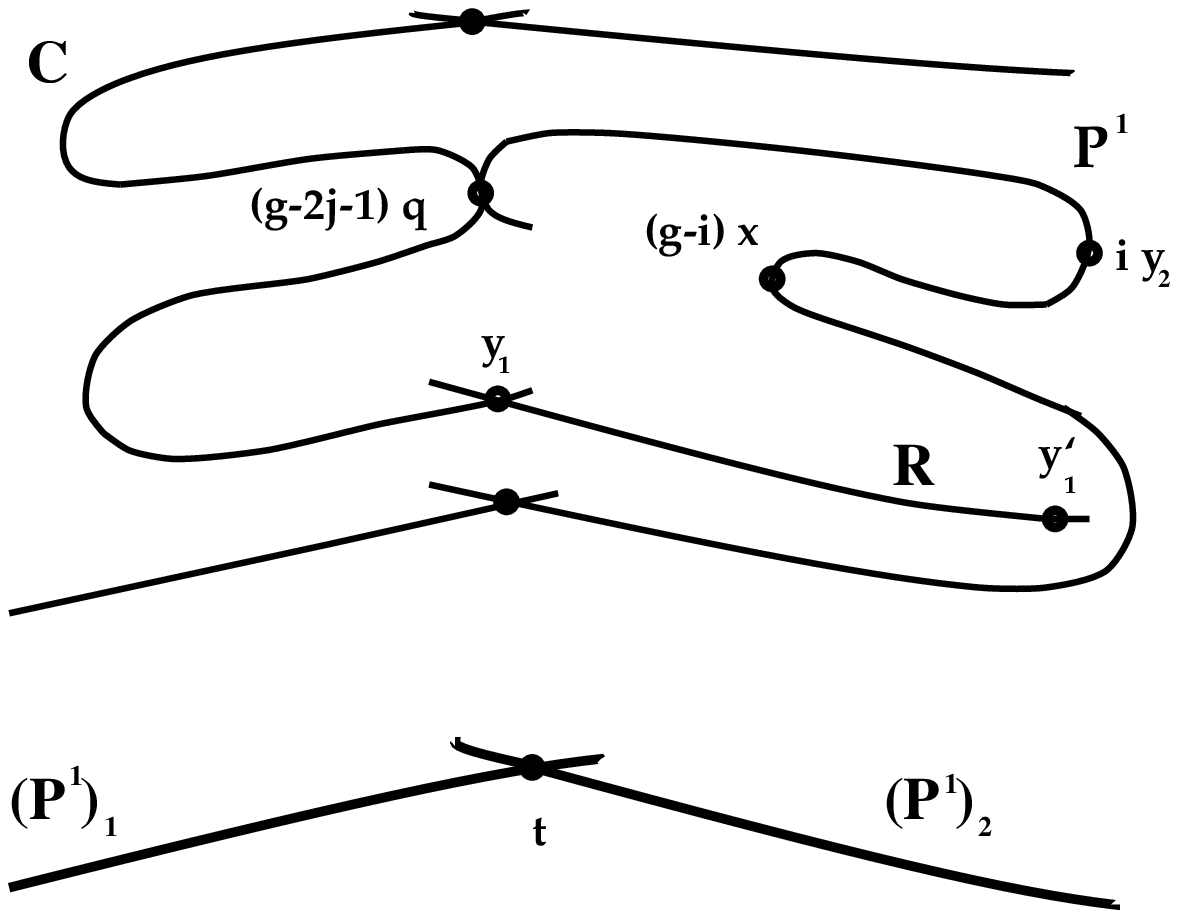,width=8cm,height=6cm,angle=0}}
\end{center}
\label{}
\caption{}
\end{figure}


\smallskip

\textbf{ii)} $w^{l_C}(q)=g-1$.
 We denote $a^{l_C}(q)=a^{l_{\PP^1}}(q)=(j,g-j)$.
 Since $\mbox{ord}_q(\tau_{\PP^1})\leq i$ we get that
 $j\leq i$. Now $w^{l_C}(q)=g-1$ is already a
 codimension $1$ condition on $\mathcal{M}_{g,1}$,
 so it follows that $\mbox{ord}_q(\sigma_C)=j$,
hence $\mbox{div}(\sigma_{\PP^1})\geq (g-j)q+iy_2$. This yields $i=j$ and
$\mbox{div}(\sigma_{\PP^1})=(g-i)q+iy_2$.
 We thus get that $[C,y_1,q]\in Y_i$. 

 Conversely, given $(C,y_1,q)\in \mathcal{M}_{g,2}$ together
 with a $\mathfrak g^1_{g-i}$ on $C$ with $(g-2i)q$ in a fiber,
 we construct a degree $g$ admissible covering
 $f:X'\rightarrow (\PP^1)_1\cup_t(\PP^1)_2$,
which will prove that $[C,q,y_2]\in (\pi_x)_*(D_{y_2}\cdot \Delta_{0:xy_2})$:
we first take $f_{|C}:C\rightarrow (\PP^1)_1$ of degree
 $g-i$ with $(g-2i)q\subseteq f_{|C}^*(t)$.
 Then $f_{|\PP^1}:\PP^1\rightarrow (\PP^1)_2$
 is of degree $g-i$, completely ramified at $x$ and with
 $f_{|\PP^1}^{-1}(t)=(g-2i)q+iy_2$. At $y_2\in \PP^1$ we insert a 
rational curve $R$ which we map $i:1$ to $(\PP^1)_1$
 such that we have total ramification both at $y_2$
 and at the point $y_2'\in R$ characterized by
$f_{|C}(y_1)=f_{|R}(y_2')$. Finally, at each of the points in
 $f_{|C}^{-1}(t)-\{q\}$ we insert a $\PP^1$
 which we map isomorphically onto $(\PP^1)_2$. 

Thus we have proved that $\mbox{supp}(\pi_x)_*(D_{y_2}\cdot \Delta_{0:xy_2})=\cup_{j=0}^i \mbox{supp}(Y_j)$. The conclusion now follows if we notice that $D_{y_2}$ is reduced and all admissible coverings we constructed are smoothable, hence $D_{y_2}\cdot \Delta_{0:xy_2}$ is reduced too.
\end{proof}

We have thus reduced the problem of computing $[D]$ to that of computing the class of all divisors $Y_j$ on $\mathcal{M}_{g,2}$ for $0\leq j\leq i$.

\begin{proposition}\label{class_Y}
For $0\leq j\leq i$ we have the following relations
 in $\rm{Pic}$$(\mathcal{M}_{g,2})_{\mathbb Q}$:
$$Y_j\equiv_{lin}a_j\lambda+b_{1j}\psi_{y_1}+b_{2j}
\psi_{y_2},\mbox{ }\mbox{   where}$$

$$a_j=-\frac{g-2j}{g}{g \choose j}+\frac{10(g-2j)}{g-2}{g-2 \choose j-1}
 \mbox{ }\mbox{ }\mbox{ for all }0\leq j\leq i,$$
$$b_{1j}=\frac{g-2j-1}{g-1}{g-1 \choose j}\mbox{ when }
j\leq i-1, b_{1i}=0, b_{2i}=\frac{(g-2i)^3-(g-2i)}{2g-2}{g \choose i}$$
 $$b_{2j}=\frac{(g-2j-1)(g^3-g^2-4g^2j+4j^2g+2jg-2j)(g-2)!}{2j!(g-1)!}
 \mbox{ for }j\leq i-1.$$
\end{proposition}

\begin{proof}
 We will compute the class of $Y_j$
 when $j\leq i-1$. The class of $Y_i$ is computed similarly.
 Let us write the following relation in $\mbox{Pic}(\overline{\mathcal{M}}_{g,2})_{\mathbb Q}$:
$$Y_j\equiv_{lin} a_j\lambda+b_{1j}\psi_{y_1}+
b_{2j}\psi_{y_2}-c_j\delta_{0:y_1y_2}+(\mbox{ other boundary terms }).$$
We use the method of test curves to determine the coefficients $a_j,b_{1j}$ and $b_{2j}$, that is, we intersect the classes appearing on both sides of the previous relation with curves inside $\overline{\mathcal{M}}_{g,2}.$ By computing intersection numbers we obtain linear relations between the coefficients $a_j,b_{1j},b_{2j}$.

By Proposition 4.5 we have that 
\begin{equation}
Z_j:=(\pi_{y_2})_*(Y_j\cdot \Delta_{0:y_1 y_2})
\equiv_{lin} a_j\lambda+c_j\psi_{y_1}+(\mbox{ boundary }).
\end{equation}
Using  the same reasoning as in Proposition~\ref{reduce},
 we obtain that $Z_j$ is the closure in $\overline{\mathcal{M}}_{g,1}$
 of the locus of curves $(C,y_1)$ carrying a
 $\mathfrak g^1_{g-j}$ with $(g-2j)y_1$ in a fiber.

In order to determine the coefficient $c_j$ we intersect both sides 
of $(5)$ with a general fiber $F$ of the map
 $\overline{\mathcal{M}}_{g,1}\rightarrow
 \overline{\mathcal{M}}_{g}$: we get that
 $c_j=Z_j\cdot F/\psi_{y_1}\cdot F=b(g-j,g)/(2g-2)$ (cf. Proposition \ref{enumerative}). 

  To determine $b_{1j}$ and $b_{2j}$ we use two
 test curves in $\overline{\mathcal{M}}_{g,2}$:
 first, we fix a general curve $C$ of genus $g$ and we obtain
 a family $C_{[1]}=\{(C,y_1,y_2)\}_{y_1\in C}$,
 by fixing a general point $y_2\in C$ and letting $y_1$ vary on $C$.
 From $(5)$, clearly
 $C_{[1]}\cdot Z_j=(2g-1)b_{1j}+b_{2j}-c_j$.
 On the other hand, according to Proposition~\ref{enumerative}
$C_{[1]}\cdot Z_j=c(g-j,g, 1)$.

 For a new relation between $b_{1j}$ and $b_{2j}$
 we use the test curve $C_{[2]}=\{(C,y_1,y_2)\}_{y_2\in C}$
 in $\overline{\mathcal{M}}_{g,2}$,
 where this time $y_1$ is a fixed general point
 while $y_2$ varies on $C$. We have the equation $(2g-1)b_
{2j}+b_{1j}-c_j=C_{[2]}\cdot Z_j=c(g-j,g,g-2j-1)$, and since
 $c_j$ is already known we get in this way both $b_{1j}$ and $b_{2j}$.

 We are only left with the computation of $a_j$.
 From \cite{EH2}, Theorem 4.1 we know that the class of
 $Z_j$ is a linear combination  of the Brill-Noether
 class and of the class of the divisor of Weierstrass points,
 that is, $Z_j\equiv_{lin} \mu
BN+\nu \mathcal{W}$, where
$$BN:=(g+3)\lambda-\frac{g+1}{6}\delta_{irr}-\sum_{i=1}^{g-1}i(g-i)
\delta_{i:y_1}\mbox{ and }$$
$$\mathcal{W}:=-\lambda+\frac{g(g+1)}{2}\psi_{y_1}-
\sum_{i=1}^{g-1}{g-i+1\choose 2}\delta_{i:y_1}.$$
We already know that $\nu=2c_j/(g(g+1))$.
 To determine $\mu$ we use the following test curve in
 $\overline{\mathcal{M}}_{g,1}$: we take a general curve $B$
 of genus $g-1$ and a general $2$-pointed elliptic curve $(E,0,y_1)$.
 We consider the family $\overline{B}=
\{X_q=B\cup_{q\sim 0}E,y_1\}_{q\in B}$ obtained by identifying
 the variable point $q\in B$ with the fixed point $0\in E$.
 We easily get $\overline{B}\cdot \psi_{y_1}=
\overline{B}\cdot \lambda=0$,
 $\overline{B}\cdot \delta_{1:y_1}=-\mbox{deg}K_B=4-2g$, 
while $\overline{B}$ vanishes on all the other boundaries.
 On the other hand $\overline{B}\cdot Z_j$ is the number of
 limit $\mathfrak g^1_{g-j}$'s on the curves $X_q$ having vanishing
 $g-2j$ at the fixed point $y_1\in E$. If $l=(l_B,l_E)$ is such a linear 
series, then using again the additivity of the Brill-Noether numbers
 (cf. \cite{EH1}, Proposition 4.6) and the assumption that $y_1-0\in
 \mbox{Pic}^0(E)$ is not torsion, we obtain that $w^{l_B}(q)=g-2j$,
 so either $a^{l_B}(q)=(1,g-2j)$ or $a^{l_B}(q)=(0,g-2j+1)$.
 Thus $\overline{B}\cdot Z_j=b(g-j-1,g-1)+b(g-j,g-1)$ and we
 can write a new relation enabling us to compute $a_j$.
\end{proof}

We can now complete the proof of Theorem~\ref{class_D}:
\begin{proof}[Proof of Theorem~\ref{class_D}]
 Let us write
 $D\equiv_{lin}A\lambda+B_1\Psi_x+B_2\Psi_y$, where
 $\Psi_x:=\sum_{j=1}^{g-i}\psi_{x_j}$ and $\Psi_{y}:=
\sum_{j=1}^{i+1}\psi_{y_j}$. As noticed before, the
 $\{\lambda, \Psi_y\}$-part of $[D]$ and the $\{\lambda, \psi_{y_1}\}$-part
 of $\sum_{j=0}^i [Y_j]$ coincide, hence using Proposition \ref{class_Y}
$$A=\sum_{j=0}^i a_j=-{g-1 \choose i}+10{g-3 \choose i-1}\mbox{ and }
 B_2=\sum_{j=0}^i b_{1j}={g-2\choose i-1}.$$ 

 Finally, to determine $B_1$ one has to compute
 the $\psi_{x_1}$ coefficient of the divisor $D_{x_2}$
 on $\overline{\mathcal{M}}_{g,3}$. Arguing in a way that
 is entirely similar to Proposition \ref{sum_Y} we obtain that
 $B_1={g-2 \choose i}$.
\end{proof}


\providecommand{\bysame}{\leavevmode\hbox to3em{\hrulefill}\thinspace}

\end{document}